\documentclass{amsart}
\usepackage{graphicx}
\theoremstyle{definition}

\theoremstyle{remark}

\numberwithin{equation}{section}

\begin{document}

\title{On the convergence exponent of the special integral
of two dimensional Tarry's problem}

\author{I.Sh.Jabbarov}
\address{AZ2000 Haydar Aliyev avenue, 159, Ganja State University, Azerbaijan}
\email{ilgar\_j@rambler.ru}

\begin{abstract}
In this article new upper and lower bounds for
the convergence exponent of the special integral
of two dimensional Tarry's problem are found.

\end{abstract}
\maketitle

\begin{center}
\noindent \textbf{1. Introduction}
\end{center}

\textbf{1. Introduction}

Let $n$ and $m$ be natural numbers, and

\[F(x,y)=\sum _{i=0}^{n}\sum _{j=0,i+j>0}^{m}\alpha _{ij} x^{i} y^{j}   \]
be a polynomial in two variables \textit{$x$ }and \textit{$y$ }with real coefficients. Under the special integral of  the two-dimensional Terry's problem one understands an integral

\[\theta _{k} =\int _{-\infty }^{\infty }\int _{-\infty }^{\infty }\cdots \int _{-\infty }^{\infty }\left|\int _{0}^{1}\int _{0}^{1} e^{2\pi iF(x,y)} dxdy \right|^{2k} d\alpha _{10} d\alpha _{01} \cdots d\alpha _{mn}    \]

\textbf{\textit{Definition.}} The number $\gamma >0$ is called to be the convergence exponent of  the integral $\theta _{k} $, if it converges for real numbers $2k<\gamma $, and diverges  for real numbers $\gamma >2k$.

The problem of finding $\gamma $ for the one-dimensional case investigated by Hua Loo Keng ([7]). This problem was completely solved in [3]. The authors introducing the notions of complete and non-complete polynomials found that the value of the convergence exponent depends on the type of the polynomial. Creating a theory of multiple trigonometric sums, the authors of the works [3 - 6] considered the multi-dimensional analogue of Tarry's problem and received the first results in the analysis of the convergence exponent for the singular integral $\theta _{k} $ of the  multi-dimensional problem.

The question on an exact value of the convergence exponent, as well as on its estimation, were considered in the works [3, 4, 8, 9, 10, 12]. In the present article we obtain new upper and lower bounds for the convergence exponent of the singular integral $\theta _{k} $.

\textbf{THEOREM 1.} The integral $\theta _{k} $ diverges for integers  \textit{$k$ }such that

\[4k\le 2+(n+m)(n+1)(m+1)/2.\]

\textbf{THEOREM 2.} The integral $\theta _{k} $ converges for integers  \textit{$k$ }such that

\[4k>2+(n+m)(n+1)(m+1)/2.\]

\textit{}
\begin{center}
\textbf{2. Auxiliary Lemmas}.
\end{center}
For the proof of our main results, we need to establish some auxiliary statements (see [10-13]).

\textbf{\textit{Lemma 1.}} Let in a bounded closed Jordan domain $\Omega $ of $n$\textit{-}dimensional space $R^{n} $ a continuous function $f(\bar{x})=f(x_{1} ,...,x_{n} )$ and continuously differentiable functions $f_{j}(\bar{x})=f(x_{1} ,...,x_{n} )$ with $j=1,...,r$ be given. Let the Jacoby matrix

\[\frac{\partial (f_{1} ,...,f_{r} )}{\partial (x_{1} ,...,x_{n} )} \]
everywhere in $\Omega $ have maximal rank. Let, further, $\bar{\xi }_{0} =(\xi _{1}^{0} ,...,\xi _{r}^{0} )$ be some interior point of the image of the mapping $\bar{x}\mapsto (f_{1} ,...,f_{r} )$  and $\bar{x}_{0} $ be any point of  $\Omega $ such that

\[f_{1} (\bar{x}_{0} )=\xi _{1}^{0} ,...,f_{r} (\bar{x}_{0} )=\xi _{r}^{0} .\]

Then, everywhere in some neighborhood of the point $\bar{\xi }_{0} $ the equality

\[\frac{\partial ^{r} }{\partial \xi _{1} \cdots \partial \xi _{r} } \int _{\Omega (\bar{\xi })}f(\bar{x}) d\bar{x}=\int _{M}f(\bar{x})\frac{ds}{\sqrt{G} }  ,\]
is satisfied; here $\Omega (\bar{\xi })$ is a sub domain in $\Omega $  defined by the system of inequalities $f_{j} (\bar{x})\le \xi _{j} $,  the surface $M(\bar{\xi })$ defined by the system of equations $f_{j} (\bar{x})=\xi _{j} \, (j=1,...,r)$, and $G$ denotes  the Gram determinant of gradients of the functions $f_{j} (\bar{x})$ , i.e. $G=\left|(\nabla f_{i} ,\nabla f_{j})\right|$.

Note that $n-r$ dimensional element of the volume on the surface of the lemma 1 we will call the element of area on the surface.

\textbf{\textit{Consequence.}} Let the conditions of Lemma 1 be satisfied. Then an equality

\[\int _{\Omega }f(\bar{x}) d\bar{x}=\int _{m_{1} }^{M_{1} }\cdots \int _{m_{r} }^{M_{r} }du_{1} \cdots du_{r} \int _{M}f(\bar{x})\frac{ds}{\sqrt{G} }    ,\]
holds, where $m_{j} $ and $M_{j} $ are, respectively, the minimal and maximal values of $f_{j} (\bar{x})$, $m_{j} \le f_{j} (\bar{x})\le M_{j} $, $M=M(\bar{u})$ is a surface in $\Omega $ defined by the system of equations $f_{j} =u_{j} ,j=1,...,r$, and $G$ denotes  the Gram determinant of gradients of the functions $f_{j} (\bar{x})$ defining the surface $M$.

\textbf{Lemma 2.} Let the polynomial $F(x,y)$ be defined by the equality of  sec. 1. Then, when $2k\ge N=(n+1)(m+1)-1$, the following formula is true

\[\theta _{k} =(2\pi )^{N} \int _{\Pi }\frac{ds}{\sqrt{G_{0} } }  ,\]
where the surface integral is taken over the surface $\Pi $, determined by the system of equations

\[x_{1} +x_{2} +\cdots +x_{k} -x_{k+1} -x_{k+2} -\cdots -x_{2k} =0,\]

\[\cdots \quad \cdots \quad \cdots \quad \cdots \]
\begin{equation} \label{1}
x_{1}^{i} y_{1}^{j} +x_{2}^{i} y_{2}^{j} +\cdots +x_{k}^{i} y_{k}^{j} -x_{k+1}^{i} y_{k+1}^{j} -x_{k+2}^{i} y_{k+2}^{j} -\cdots -x_{2k}^{i} y_{2k}^{j} =0,
\end{equation}

\[\cdots \quad \cdots \quad \cdots \quad \cdots \]

\[x_{1}^{n} y_{1}^{m} +x_{2}^{n} y_{2}^{m} +\cdots +x_{k}^{n} y_{k}^{m} -x_{k+1}^{n} y_{k+1}^{m} -x_{k+2}^{n} y_{k+2}^{m} -\cdots -x_{2k}^{n} y_{2k}^{m} =0\]
in the-$4k$dimensional unite cube;

1) here $N=(n+1)(m+1)-1$ is the number of  monomials of the polynomial $F(x,y)$, and $G_{0} $ means the Gram determinant of gradients of  functions standing on the left parts of the system (1), i.e. $G_{0} =\det (A_{0} \cdot {}^{t} A_{0} )$ and

\[A_{0} =\left(\begin{array}{ccccc} {1} & {0} & {\cdots } & {-1} & {0} \\ {0} & {1} &  {\cdots } & {0} & {-1} \\ {\cdots } & {\cdots } &  {\cdots } & {\cdots } & {\cdots } \\ {ix_{1}^{i-1} y_{1}^{j} } & {jx_{1}^{i} y_{1}^{j-1} } &  {\cdots } & {-ix_{2k}^{i-1} y_{2k}^{j} } & {-jx_{2k}^{i} y_{2k}^{j-1} } \\ {\cdots } & {\cdots } &  {\cdots } & {\cdots } & {\cdots } \\ {nx_{1}^{n-1} y_{1}^{m} } & {mx_{1}^{n} y_{1}^{m-1} } &  {\cdots } & {-nx_{2k}^{n-1} y_{2k}^{m} } & {-mx_{2k}^{n} y_{2k}^{m-1} } \end{array}\right);\]

2) the equality of the lemma is understood in the meaning that both sides of this equality converges and diverges simultaneously, and convergence of the surface integral is defined as follows:
\[(2\pi )^{N} \int _{\Pi }\frac{ds}{\sqrt{G_{0} } } =lim_{\eta\rightarrow 0} (2\pi )^{N} \int _{\Pi, \bar{x}\in D_{\eta}\times D_{\eta} }\frac{ds}{\sqrt{G_{0} } }\]
here $D_{\lambda } \subset [0,1]^{2k} $ is a subdomain defined by the equality $D_{\lambda } =\left\{\left(\bar{x}_{1} ,...,\bar{x}_{k} \right)|G\ge \lambda \right\}$, \textit{G} is a Gram determinant of gradients of the functions
\[f_{ij} =x_{1}^{i} y_{1}^{j} +x_{2}^{i} y_{2}^{j} +\cdots +x_{k}^{i} y_{k}^{j}.\]

We need in the following multidimensional analogue of the theorem of Arzela (see. [43, p. 752]), which can easily be proven by the method of mathematical induction.

\textbf{\textit{Lemma 3.}} Let we are given with the sequence of functions

\[f_{n} (\bar{x})\, (n=1,2,...),\]
integrable in the product $K=[a_{1} ,b_{1} ]\times [a_{2} ,b_{2} ]\times \cdots \times [a_{s} ,b_{s} ]$ and bounded in their totality

\[\left|f_{n} (\bar{x})\right|\le L(\bar{x}\in K,n=1,2,...).\]
Let for all $\bar{x}\in K$\textit{ }there is a limit

\[\varphi (\bar{x})=\mathop{\lim }\limits_{n\to \infty } f_{n} (\bar{x}).\]
If for any $r,0\le r\le s-1$ functions $f_{n} (\bar{x})$ and $\varphi (\bar{x})$\textit{ }are integrable on

\[[a_{r+1} ,b_{r+1} ]\times \cdots \times [a_{s} ,b_{s} ],\]
then

\[\mathop{\lim }\limits_{n\to \infty } \int _{K}f_{n} (\bar{x})d\bar{x} =\int _{K}\varphi (\bar{x}) d\bar{x}.\]

We now introduce a matrix  $A_{1} $  which we get by arranging of the entries of columns of the matrix $A_{0} $, consequently, in a line, with subsequently taking of the transposed Jacoby matrix of obtained system of functions. Analogically, we introduce matrices  $A_{2} ,...,A_{m+n-1} $.

Let's consider a polynomial given by the equality
\begin{equation}\label{2}
F(x,y)=\sum _{j=1}^{N}\alpha _{j} \gamma _{j} (x,y) ,
\end{equation}
where $\gamma _{j} (\bar{x})=\gamma _{j} (x,y)$ are the monomials of a kind

\[\gamma _{j} (\bar{x})=x^{k_{j} } y^{l_{j} } .\]

The system connected with this polynomial defines a mapping $\varphi ={}^{t} \left(\varphi _{1} ,\varphi _{2} ,...,\varphi _{N} \right)$, where

\[\varphi _{j} :(x,y)\mapsto \left(\begin{array}{ccc} {\gamma _{j} (x_{1} ,y_{1} )+} & {\cdots } & {-\gamma _{j} (x_{2k} ,y_{2k} )} \end{array}\right).\]

 The Jacoby matrix looks like:

\[\left(\begin{array}{cccc} {k_{1} x_{1}^{-1} \gamma _{1} (\bar{x}_{1} )} & {l_{1} y_{1}^{-1} \gamma _{1} (\bar{x}_{2} )} & {\cdots } & {-l_{1} y_{2k}^{-1} \gamma _{1} (\bar{x}_{2k} )} \\ {k_{2} x_{1}^{-1} \gamma _{2} (\bar{x}_{1} )} & {l_{2} y_{1}^{-1} \gamma _{2} (\bar{x}_{2} )} & {\cdots } & {-l_{2} y_{2k}^{-1} \gamma _{2} (\bar{x}_{2k} )} \\ {\vdots } & {\vdots } & {\ddots } & {\vdots } \\ {k_{N} x_{1}^{-1} \gamma _{N} (\bar{x}_{1} )} & {l_{N} y_{1}^{-1} \gamma _{N} (\bar{x}_{2} )} & {\cdots } & {-l_{N} y_{2k}^{-1} \gamma _{N} (\bar{x}_{2k} )} \end{array}\right).\]

Minor of order 2 of this matrix, composed of the first 2 columns, after of reducing by common factors of the elements of columns and lines of the determinant, lead to minors of a matrix of a following kind
\begin{equation}\label{3}
\left(\begin{array}{cc} {\begin{array}{c} {k_{1} } \\ {k_{2} } \\ {\vdots } \\ {k_{N} } \end{array}} & {\begin{array}{c} {l_{1} } \\ {l_{2} } \\ {\vdots } \\ {l_{N} } \end{array}} \end{array}\right).
\end{equation}
Singularities of the mapping $\varphi =\left(\varphi _{1} ,\varphi _{2} ,...,\varphi _{N} \right)$ depend on the rank of this matrix. If the rank is maximal then the given mapping is regular everywhere, with exception of the points of a set of  zero measure for $k\ge 1$.

\textbf{\textit{ Lemma 4. }}Let the matrix (3) have the rank $\rho \ge 1$ and let its first $\rho $ columns be linearly independent. We will designate

\[K(\bar{x})=\left(\begin{array}{cc} {\begin{array}{c} {k_{1} \gamma _{1} } \\ {k_{2} \gamma _{2} } \\ {\vdots } \\ {k_{N} \gamma _{N} } \end{array}} & {\begin{array}{c} {l_{1} \gamma _{1} } \\ {l_{2} \gamma _{21} } \\ {\vdots } \\ {l_{N} \gamma _{N} } \end{array}} \end{array}\right).\]

Then, for all natural $k$ such that $k\rho \ge N$ the block matrix

\[\left(K(\bar{x}_{1} )\, \, K(\bar{x}_{2} )\, \, \cdots \, \, K(\bar{x}_{k} )\right)\]
has a rank $N$ for all $\left(\bar{x}_{1} ,\bar{x}_{2} ,\cdots ,\bar{x}_{k} \right)$, with exception for the points belonging to a subset of zero Jourdan measure.

\textit{  Proof. }It is visible immediately that for the values of variables distinct from zero the matrix $K(\bar{x})$ obtained, beginning from the matrix (3), by applying elementary transformations over columns and lines has a rank $\rho $ also. In any open set of varying of a vector $\bar{x}_{2} $ the columns of the matrix $K(\bar{x}_{2} )$ cannot belong to the linear span of columns of the matrix $K(\bar{x}_{1} )$. Really, let for a given $\bar{x}_{2} $ the first column of the matrix $K(\bar{x}_{2} )$ be a linear combination of the first $\rho $ linear independent columns of the matrix $K(\bar{x}_{1} )$. Then, the coefficients of this linear combination (depending on $\bar{x}_{2} $) satisfy the system of equations:
\begin{equation}\label{4}
K(\bar{x}_{1} )\bar{\lambda }(\bar{x}_{2} )=\left(K(\bar{x}_{2} )\right)_{1} ,
\end{equation}
where at the right hand side of the equality the first column of the matrix $K(\bar{x}_{2} )$ stands, and$\bar{\lambda }(\bar{x}_{2} )\in R^{\rho } $ denotes some vector. Designating $K$ a base minor of the matrix $K(\bar{x}_{1} )$ which is assumed to be located at the left top corner, and $\bar{k}(\bar{x}_{2} )$ is a vector with coordinates coincident with the first $\rho $ components of the first column of the matrix $K(\bar{x}_{2} )$, we will have:

\[\bar{\lambda }(\bar{x}_{2} )=K^{-1} \bar{k}(\bar{x}_{2} );\]
here we simultaneously designate by $K$ a matrix of a base minor also. From this relation it is visible that the coordinates of the vector $\bar{\lambda }(\bar{x}_{2} )$ are rational functions of $\bar{x}_{2} $. Substituting the found vector into (4), we get:

\[K(\bar{x}_{1} )\{ K^{-1} \bar{k}(\bar{x}_{2} )\} =\left(K(\bar{x}_{2} )\right)_{1} .\]
It means that all elements of the first column of the matrix $K(\bar{x}_{2} )$, for example, the monomial $k_{N} \gamma _{N} (\bar{x}_{2} )$ is expressed by the monomials $k_{1} \gamma _{1} (\bar{x}_{2} )$, $k_{2} \gamma _{2} (\bar{x}_{2} )$ linearly with coefficients being rational functions of $\bar{x}_{1} $ in an open set of changing $\bar{x}_{2} $, which is impossible. As the algebraic set is a Jourdan set, then the points where the first column of the matrix $K(\bar{x}_{2} )$ belongs to a linear span of columns of the matrix $K(\bar{x}_{1} )$, by the told above, consists of boundary points, and consequently has a zero Jourdan measure.

Further, as the elements of columns of the matrix $K(\bar{x}_{2} )$ differ each from other by  numerical multipliers only, then from here it follows that any non-trivial linear combination of columns of the matrix $K(\bar{x}_{2} )$, with coefficients don't dependent on $\bar{x}_{2} $, cannot belong to a linear span of columns of the matrix $K(\bar{x}_{1} )$ in an open set of changing $\bar{x}_{2} $. So, the matrix $\left(K(\bar{x}_{1} )\, \, (K(\bar{x}_{2} ))_{1} \right)$, which is constructed by joining to the entries of the matrix $K(\bar{x}_{1} )$ the first column of the matrix $K(\bar{x}_{2} )$ from the right hand side,  has the maximal rank everywhere, except for points of a Jourdan set of zero measure. We will prove that the same statement is true for the matrix $\left(K(\bar{x}_{1} )\, \, (K(\bar{x}_{2} ))_{1} \, \, (K(\bar{x}_{2} ))_{2} \right)$ when the matrix (3) has a rank 2. Let us suppose, in the contrary, that the last column of this matrix belongs to the linear span of the previous columns. Arguing as above, we find analogical relation to the equality (4):
\begin{equation}\label{5}
K_{2} \bar{\lambda }_{2} (\bar{x}_{2} )=\left(K(\bar{x}_{2} )\right)_{2} ,
\end{equation}
where $K_{2} $ is formed of lines on which the base minor of the matrix $\left(K(\bar{x}_{1} )\, \, (K(\bar{x}_{2} ))_{1} \right)$ is located, and $\bar{\lambda }_{2} (\bar{x}_{2} )$ is a vector made of coefficients of a linear combination by means of which the last column is expressed linearly by previous vectors. Taking the first 3 components (it is specified by the top index 3) of the columns on both parts of the equality (5) we receive:

\[K_{2}^{3} \bar{\lambda }_{2} (\bar{x}_{2} )=\left(K(\bar{x}_{2} )\right)_{2}^{3} .\]

Let's consider the last of the components of the vector ${}^{t} \bar{\lambda }_{2} (\bar{x}_{2} )=(\lambda _{1} ,...,\lambda _{3} )$. By Kramer's rule, it is represented as a relation $\lambda _{3} =D'/D$ where $D$ is a base minor, and $D'$ is gotten from this minor  replacing the last column by the column $\left(K(\bar{x}_{2} )\right)_{2}^{3} $. Then, as above, we receive that any monomial, say, the monomial $l_{N} \gamma _{N} (\bar{x}_{2} )$ is expressed as a linear combination of monomials $l_{1} \gamma _{1} (\bar{x}_{1} )$, $l_{2} \gamma _{1} (\bar{x}_{2} )$, $l_{3} \gamma _{3} (\bar{x}_{2} )$. We have:

\[l_{N} \gamma _{N} (\bar{x}_{2}) )=\lambda _{1} k_{N} \gamma _{N} (\bar{x}_{1} )+\lambda _{2} l_{N} \gamma _{N} (\bar{x}_{1} )+\lambda _{3} k_{N} \gamma _{N} (\bar{x}_{2} ),\]
or
\begin{equation}\label{6}
(l_{N} -\lambda _{3} k_{N} )\gamma _{N} \left(\bar{x}_{2} \right )=\lambda _{1} k_{N} \gamma _{N} \left (\bar{x_{1}} \right )+\lambda _{2} l_{N} \gamma _{N} \left(\bar{x}_{1} \right ).       \end{equation}

Obviously that

\[D=\left|\begin{array}{ccc} {k_{1} \gamma _{1} (x_{1} ,y_{1} )} & {l_{1} \gamma _{1} (x_{1} ,y_{1} )} & {k_{1} \gamma _{1} (x_{2} ,y_{2} )} \\ {k_{2} \gamma _{2} (x_{1} ,y_{1} )} & {l_{2} \gamma _{2} (x_{1} ,y_{1} )} & {k_{2} \gamma _{2} (x_{2} ,y_{2} )} \\ {k_{3} \gamma _{3} (x_{1} ,y_{1} )} & {l_{3} \gamma _{3} (x_{1} ,y_{1} )} & {k_{3} \gamma _{3} (x_{2} ,y_{2} )} \end{array}\right|\]

and

\[l_{N} -\lambda _{3} k_{N} =\frac{l_{N} D-k_{N} D'}{D} .\]

Designating columns of the determinant $D$ as $D_{1} $, \dots , $D_{3} $ we can write the nominator of the last fraction as follows:

\[l_{N} D-k_{N} D'=\det (D_{1} ,D_{2} ,l_{N} \left(K(\bar{x}_{2} )\right)_{1}^{3} -k_{N} \left(K(\bar{x}_{2} )\right)_{2}^{3} ).\]

In the consent with the told above, the last column cannot belong to the linear span of the previous columns. Therefore, the coefficient of the monomial on the left part of (6) doesn't vanish in the same open set everywhere, with exception for points of a subset of zero Jourdan measure. Further,

\[l_{N} \left(K(\bar{x}_{2} )\right)_{1}^{3} -k_{N} \left(K(\bar{x}_{2} )\right)_{2}^{3} =\left(\begin{array}{c} {(l_{N} k_{1} -k_{N} l_{1} )\gamma _{1} (x_{2} ,y_{2} )} \\ {(l_{N} k_{2} -k_{N} l_{2} )\gamma _{2} (x_{2} ,y_{2} )} \\ {(l_{N} k_{3} -k_{N} l_{3} )\gamma _{3} (x_{2} ,y_{2} )} \end{array}\right),\]
and all of the coefficients of monomials aren't zero. Rewriting the equality (6) in the view
\begin{equation}\label{7}
(l_{N} D-k_{N} D')\gamma _{N} (x_{2} ,y_{2} )=D(\lambda _{1} k_{N} +\lambda _{2} l_{N} )\gamma _{N} (\bar{x}{}_{1} ),
\end{equation}
we note that the coefficient at the monomial $\gamma _{N} (\bar{x}_{2} )$ on the left hand side of the equality (7) is a polynomial of monomials $\gamma _{1} (\bar{x}_{2} )$,$\gamma _{2} (\bar{x}_{2} )$,$\gamma _{3} (\bar{x}_{2} )$ (with coefficients which are polynomials of the variable $\bar{x}_{1} $), and is distinct from zero in any open set everywhere, with exception for the points of a subset of zero Jourdan measure.  From the equality it is obvious that the expression standing on the right hand side of the (7) such a polynomial is also (with another coefficients). Then, taking as a monomial $\gamma _{N} (\bar{x}_{2} )$ that one of monomials which has a higher degree, we see that the left and right hand sides of the (6) are the polynomials of different degree of the variable $\bar{x}_{2} =(x_{2} ,y_{2} )$. So, the equality (7) changes now into the polynomial equation that has a set of solution having zero Jourdan measure.

So, the matrix  $\left(K(\bar{x}_{1} )\, \, K(\bar{x}_{2} )\right)$ has a maximal rank for all pares $\left(\bar{x}_{1} ,\bar{x}_{2} \right)$ with exception for the points of a set of zero Jourdan measure. Joining blocks consequently, we arrive at the proof of the lemma 4.

 \textbf{\textit{Lemma 5. }} All of matrices $A_{j+1} $ ($j\ge 0$)  have maximal rank for all values of variables, with exception for their values from a subset of zero Jourdan measure.

 As it is visible from the definition of the matrices $A_{j+1} $ ($j\ge 0$), entries of these matrices constructed of the blocks of order 2 which are transposed Jacoby matrices of gradient of monomials of a view $x^{i} y^{j} \, (i+j>1)$ (with exception of a case $i+j=m+n-1$). So, such a block has a determinant:

\[\left|\begin{array}{cc} {i(i-1)x^{i-2} y^{j} } & {ijx^{i-1} y^{j-1} } \\ {ijx^{i-1} y^{j-1} } & {j(j-1)x^{i} y^{j-2} } \end{array}\right|=ij(1-i-j)x^{2i-2} y^{2j-2} ,\]
which is non-zero when both of variables doesn't vanishing. By this reason all of matrices $A_{j+1} $ have maximal rank for all values of variables, with exception for values from subset of zero Jourdan measure.

\textit{Proof of Lemma 2.} To the integral

\[\left(\int _{0}^{1}\int _{0}^{1}e^{2\pi iF(x,y)} dxdy  \right)^{k} =\]
\[\int _{0}^{1}\cdots \int _{0}^{1}e^{2\pi i\sum _{i=0}^{n}\sum _{j=0,i+j>0}^{m}\alpha _{ij} ( x_{1}^{{\rm i}} y_{1}^{j} +\cdots +x_{k}^{{\rm i}} y_{k}^{j} ) } dx_{1} \cdots dy_{k}   \]
apply the corollary to Lemma 1, taking as the functions $f_{j} $ the polynomials

\[f_{ij} =x_{1}^{i} y_{1}^{j} +\cdots +x_{k}^{i} y_{k}^{j} .\]

In consent with the lemma 4 the equality $G=0$ can only be satisfied on a subset of smaller dimension in the $2k$-dimensional unite cube and, therefore, has a zero measure. Applying outside of this manifold in every compact $D_{\eta}\times D_{\eta} $ the corollary to Lemma 1, we obtain the following equality:
\[\int _{0}^{1}\cdots \int _{0}^{1}e^{2\pi i\sum _{i=0}^{n}\sum _{j=0,i+i>0}^{m}\alpha _{ij} (x_{1}^{i} y_{1j}^{} +\cdots x_{k}^{i} y_{k}^{j} )  } dx_{1} \cdots dy_{k} = \]
\begin{equation}\label{8}
\int _{0}^{k}\cdots \int _{0}^{k}\left(\int _{\Pi (\bar{u})}\frac{ds}{\sqrt{G} }  \right)   e^{2\pi i\sum _{i=0}^{n}\sum _{j=0,i+j>0}^{m}\alpha _{ij} u_{ij}   } du_{10} \cdots du_{nm};
\end{equation}
here $\Pi (\bar{u})$ is a surface defined by the system of equations $f_{ij} =u_{ij} ;0\le i\le n$ $0\le j\le m,$ $i+j>0.$ Then, considering the last integral as a Fourier transformation, we have by virtue of  Parseval equality:
\begin{equation}\label{9}
\theta _{k} =(2\pi )^{N} \int _{0}^{k}\cdots \int _{0}^{k}\left(\int _{\Pi (\bar{u})}\frac{ds}{\sqrt{G} }  \right)  ^{2} du_{10} \cdots du_{nm} , \end{equation}
and the equality is understood in the sense that if one of  two parts of the last equality converges then the other part converges also, and the corresponding values are equal.

Now, let us turn to the question of exact justification of equality (9). We prove that in (9) the integral on the right part exists as improper integral in the cube $[0,k]^{N} $. Consider the mapping $\bar{f}:[0,1]^{2k} \to [0,k]^{N} $, where $\bar{f}=\left(f_{ij} \right)$ is a vector - function components of  which coincide with polynomials defining $\Pi(\bar{u})$. Let's consider some decreasing sequence of positive numbers $\left(\eta _{l} \right)$ ($\eta _{l} \to 0$). The union of closed sets $U_{l} =\bar{f}\left(\{ \bar{x}\in [0,1]^{2k} |G\ge \eta _{l} \} \right)$, ($l=1,2,...$) contains $[0,k]^{N} \backslash \bar{f}\left(W\right)$, and their prototypes $\bar{f}^{-1} \left(U_{l} \right)$ contains the  region $G>0$. In each of  closed region $G\ge \eta _{l} >0$ apply corollary to Lemma 1:

\[\int _{0}^{1}\cdots \int _{0,G\ge \eta _{l} }^{1}e^{2\pi i\sum _{i=0}^{n}\sum _{j=o,i+i>0}^{m}\alpha _{ij} (x_{1}^{i} y_{1}^{j} +\cdots x_{k}^{i} y_{k}^{j} )  } dx_{1} \cdots dy_{k} = \]
\[\int _{0}^{k}\cdots \int _{0}^{k}\left(\int _{\Pi (\bar{u}),\bar{u}\in U_{l} }\frac{ds}{\sqrt{G} }  \right)   e^{2\pi i\sum _{i=0}^{n}\sum _{j=0,i+j>0}^{m}\alpha _{ij} u_{ij}   } du_{10} \cdots du_{nm} .\]

Since the left hand side, when $\eta _{l} \to 0$, tends to the value of the integral on the left part of the formula (2), we obtain (2) in improper sense. Thus,  (3) takes place in improper  sense also.

Now we will prove that (3) has Lemma 2 as a consequence. Let $D$ be a closed subdomain of the unite cube where $G>0$ . Then, by Lemma 1, we have:

\begin{equation} \label{10} \int _{\bar{x}\in D,\Pi (\bar{u})}\frac{ds}{\sqrt{G} }  =\mathop{\lim }\limits_{h\to 0} \frac{1}{(2h)^{N} } \int _{u_{ij} -h<f_{ij} <u_{ij} +h}d\bar{x} ;\bar{x}=(x_{1} ,...,y_{k} ). \end{equation}

Dissect the region \textit{D}, as above, into closed parts $D_{\nu } ,\nu =1,...,l$, in each of which some minor has a maximal value for the modulus among others, and shall make the change of variables the Jacobean of which has the inverse values being equal to the maximal minor. Then, the surface integral splits into a sum of integrals of  bounded functions. Therefore, designating by $\varphi _{D} (\bar{u})$ the left side of (4) , we may, by Lemma 1 and the corollary to it, write

\[\int _{0}^{k}\cdots \int _{0}^{k}\varphi _{D} (\bar{u})\varphi _{D} (\bar{u})  d\bar{u}=\]
\[\int _{0}^{k}\cdots \int _{0}^{k}\varphi _{D} (\bar{u})\left(\mathop{\lim }\limits_{h\to 0} \frac{1}{(2h)^{N} } \sum _{\nu } \int _{u_{ij} -h\le f_{ij} \le u_{ij} +h,\bar{x}\in D_{\nu } }d\bar{x} \right)  d\bar{u}.\]

According to Lemma 3 under the integral on the right part, we can rearrange the order of integration and pass to the limit. To do this, one needs to put $h=h_{n} $ where $h_{n} \to 0$, and apply Lemma 3 to our integral at $h=h_{n} $:

\[\int _{0}^{k}\cdots \int _{0}^{k}\varphi _{D} (\bar{u})\left(\mathop{\lim }\limits_{n\to \infty } \frac{1}{(2h_{n} )^{N} } \sum _{\nu }\int _{u_{ij} -h_{n} \le f_{ij} \le u_{ij} +h_{n} ,\bar{x}\in D_{\nu } }d\bar{x}  \right)  d\bar{u}.\]

Then we have:

\[\int _{0}^{k}\cdots \int _{0}^{k}\varphi _{D} (\bar{u})\left(\mathop{\lim }\limits_{h\to 0} \frac{1}{(2h)^{N} } \int _{u_{ij} -h\le f_{ij} \le u_{ij} +h,\bar{x}\in D_{\nu } }d\bar{x} \right)  d\bar{u}=\]

\[\mathop{\lim }\limits_{h\to 0} \frac{1}{(2h)^{N} } \int _{0}^{k}\cdots \int _{0}^{k}\varphi _{D} (\bar{u})\left(\int _{u_{ij} -h\le f_{ij} \le u_{ij} +h,\bar{x}\in D_{\nu } }d\bar{x} \right)  d\bar{u}.\]

Summing up over $\nu $, we have:

\[\int _{0}^{k}\cdots \int _{0}^{k}\varphi _{D} (\bar{u})\varphi _{D} (\bar{u})  d\bar{u}=\]
\begin{equation}\label{11}
\int _{0}^{k}\cdots \int _{0}^{k}\varphi _{D} (\bar{u})\left(\mathop{\lim }\limits_{h\to 0} \frac{1}{(2h)^{N} } \int _{u_{ij} -h\le f_{ij} \le u_{ij} +h}d\bar{x} \right)  d\bar{u}=
\end{equation}
\[=\mathop{\lim }\limits_{h\to 0} \frac{1}{(2h)^{N} } \int _{0}^{k}\cdots \int _{0}^{k}\int _{\bar{x}'\in D,\Pi '(\bar{u})}\left(\int _{u_{ij} -h\le f_{ij} \le u_{ij} +h}d\bar{x} \right)\frac{ds'}{\sqrt{G'} }    d\bar{u};\]
here $ds'$ is an element of  area on the surface $\Pi '(\bar{u})$, defined in $D$ by the system of equations, $f_{ij} (\bar{x}')=u_{ij} ;$ $\; \; 0\le i\le n,\; 0\le j\le m,$ $i+j>0,$ and $G'$ has a similar meaning. Consider, for a fixed \textit{h} , the inner integral in the last chain of equalities (5), i.e., the integral

\[\int _{\bar{x}'\in D,\Pi '(\bar{u})}\left(\int _{u_{ij} -h\le f_{ij} \le u_{ij} +h}d\bar{x} \right)\frac{ds'}{\sqrt{G'} }  .\]

For the points $\bar{x}'\in [0,1]^{2k} $ we define a function $f(\bar{x}')$, taking its value  when $\bar{x}'\in \Pi '\left(\bar{u}\right)$ equal to the inner integral. We prove that the function $f(\bar{x}')$, defined by this way is continuous on $D$. Let $\bar{x'}_{1} ,\bar{x'}_{2} \in D$, $\bar{x}'_{1} =(x'_{11} ,...,y'_{1k} ),\, \bar{x}'_{2} =(x'_{21} ,...,y'_{2k} ),\, $ $\sum [(x'_{1i} -x'_{2i} )^{2} +(y'_{1i} -y'_{2i} )^{2} ] \le \varepsilon $, if $\varepsilon >0$ is specified. Then, assuming $u_{ij}^{1} =f_{ij} (\bar{x}'_{1} ),\, \quad u_{ij}^{2} =f_{ij} (\bar{x}'_{2} )$ (here we use the upper indexing), according to the mean value theorem, for some $\bar{\theta }$ and $\bar{\eta }$ we get:

\[\left|u_{ij}^{1} -u_{ij}^{2} \right|=\left|\sum _{r=1}^{k}\left(\frac{\partial f_{ij} }{\partial x_{r} } (\bar{x}_{1} ^{{'} } +\bar{\theta })\left(x'_{1r} -x'_{2r} \right)+\frac{\partial f_{ij} }{\partial y_{r} } (\bar{y}_{1} ^{{'} } +\bar{\eta })\left(y'_{1r} -y'_{2r} \right)\right) \right|\le \]

\[\le (m+n)^{3/2} (\varepsilon )^{1/2} ,\]
if $\sum [(x'_{1i} -x'_{2i} )^{2} +(y'_{1i} -y'_{2i} )^{2} ] \le \varepsilon $. Therefore, we find:

\[\left|f(\bar{x}'_{1} )-f(\bar{x}'_{2} )\right|\le \]
\[\le 2^{N} \mathop{\max }\limits_{i,j} \left\{\int _{u_{ij}^{1} -\sqrt{2k\varepsilon } \le f_{ij} \le u_{ij}^{2} +\sqrt{2k\varepsilon } }d\bar{x} +\int _{u_{ij}^{1} +h-\sqrt{2k\varepsilon } \le f_{ij} \le u_{ij}^{2} +h+\sqrt{2k\varepsilon } }d\bar{x} \right\}.\]

One estimates these integrals equally. Denoting modulus of maximal minor  as $\left|J\right|$, we get $G_{1} \le G\le C_{2r}^{N} \left|J\right|^{2} $, where $G_{1} $ is the minimal value of \textit{G} in \textit{D}. Because the independent variables are changed in $[0,1]$, using the expression for the element of area was found in the proof of Lemma 1, we find:

\[\int _{u_{ij}^{1} -\sqrt{2k\varepsilon } \le f_{ij} \le u_{ij}^{2} +\sqrt{2k\varepsilon } }d\bar{x} =\int _{u_{ij}^{1} -\sqrt{2k\varepsilon } \le z_{ij} \le u_{ij}^{2} +\sqrt{2k\varepsilon } }dz_{01} \cdots dz_{nm} \int _{\Pi \left(\bar{z}\right)}\frac{ds}{\sqrt{G} } \le   \]

\[\le C_{2k}^{N} G_{1}^{-1/2} (8\varepsilon )^{1/2} .\]

Counting arbitrariness of $\varepsilon >0$, from this we conclude continuity of $f(\bar{x}')$. Thus,

\[\int _{\Pi '(\bar{u})}\left(\int _{u_{ij} -h\le f_{ij} \le u_{ij} +h}d\bar{x} \right)\frac{ds'}{\sqrt{G'} }  =\int _{\Pi '\left(\bar{u}\right)}f\left(\bar{x}'\right)\frac{ds'}{\sqrt{G'} }  .\]

Applying corollary to Lemma 1, we find out:

\[\int _{0}^{k}\cdots \int _{0}^{k}\int _{\bar{x}'\in D,\Pi '(\bar{u})}f(\bar{x}')\frac{ds'}{\sqrt{G'} }    d\bar{u}=\int _{0}^{1}\cdots \int _{0}^{1}f(\bar{x}')d\bar{x}'=  \]

\[=\int _{0}^{1}dx'_{1} \cdots \int _{0}^{1}dy'_{k} \int _{u_{ij} -h\le x_{1}^{i} y_{1}^{j} +\cdots +x_{k}^{i} y_{k}^{j} \le u_{ij} +h}dx_{1} \cdots dy_{k} .   \]

Taking into account that $u_{ij} =x\prime_{1}^{i} y\prime_{1}^{j} +\cdots +x\prime_{k}^{i} y\prime_{k}^{j} $, in new notations with $x_{k+s} =x'_{s} ,y_{k+s} =y'_{s} $, from (5) we obtain finally:

\[\int _{0}^{k}\cdots \int _{0}^{k}\varphi _{D} (\bar{u})  \varphi _{D} (\bar{u})d\bar{u}=\]

\[=\mathop{\lim }\limits_{h\to 0} \frac{1}{(2h)^{N} } \int _{u_{ij} -h\le x_{1}^{i} y_{1}^{j} +\cdots +x_{k}^{i} y_{k}^{j} \le u_{ij} +h}dx_{1} \cdots dy_{k} =\]
\[\int _{\bar{x}\in D\times D,\Pi _{0} }\frac{ds}{\sqrt{G_{0} } } ;\bar{x}=(x_{1} ,...,y_{2k} ) , \]
where $\Pi _{0} $ is a manifold of  solutions of the system (1), and $G_{0} $ is a Gram determinant of the gradients of the functions $f_{ij} -f'_{ij} $, standing on the left parts of the system (1). From the expansion for the Gram determinant (see. [17, p. 245]), it follows that $G_{0} \ge G+G'>0$ in $D\times D$. 

It is easy to observe that the all of  reasonings of the proof above can be inverted in opposite direction. Lemma 2 is proved.

\textbf{Lemma 6}. Under the conditions of Lemma 2 the formula

\[\theta _{k} =2^{N} \Gamma (1+N/2)\pi ^{N/2} \int _{-\infty }^{\infty }\cdots \int _{-\infty }^{\infty }\left(\int _{\Pi (\bar{\alpha })}ds \right) d\alpha _{01} \cdots d\alpha _{nm} , \]
holds when  $\bar{\alpha }=(\alpha _{01} ,...,\alpha _{nm} ),\; $and $\Pi \left(\bar{\alpha }\right)$ is a part of the surface defined in the 4k- dimensional unite cube by the system of equations (1) and an inequality

\[\sum _{j=1}^{2k}\left(\left|\frac{\partial F(x_{j} ,y_{j} )}{\partial x_{j} } \right|^{2} +\left|\frac{\partial F(x_{j} ,y_{j} )}{\partial y_{j} } \right|^{2} \right) \le 1;\]
here $\Gamma $ is the Euler's gamma function, and this formula is understood in the sense that the both sides of the equality converges, if one of two parts of the equality is convergent.

\textit{Proof.} From the known relationship of [13, p. 131, the problem 35] it follows that if as a matrix \textit{A} take the matrix $A_{0} \cdot {}^{t} A_{0} $ then

\[\left(\sqrt{G_{0} } \right)^{-1} =\frac{\Gamma (N/2+1)}{\pi ^{N/2} } \int _{\left\| {}^{t} A\bar{\alpha }\right\| \le 1}d\alpha _{01} \cdots d\alpha _{nm} , \]
where $A_{0} $ means the Jacoby matrix of the system of functions $f_{ij} $. Substituting this expression in the Lemma 2 for the value of the Gram determinant, and using nonnegativity of the functions under the sighn of integral, will have:

\[\theta _{k} =(2\pi )^{N} \int _{\Pi }\frac{ds}{\sqrt{G_{0} } }  =2^{N} \pi ^{N/2} \Gamma (N/2+1)\int _{\Pi }ds\int _{\left\| {}^{t} A\bar{\alpha }\right\| \le 1}d\alpha _{01} \cdots d\alpha _{nm}   .\]

Change the order of integrations in the right part where the surface integral is taken over that part of the surface $\Pi $ on which the  inequality

\[\left(\bar{\alpha },A_{0} \cdot {}^{t} A_{0} \bar{\alpha }\right)=\sum _{j=1}^{2k}\left(\left|\frac{\partial F(x_{j} ,y_{j} )}{\partial x_{j} } \right|^{2} +\left|\frac{\partial F(x_{j} ,y_{j} )}{\partial y_{j} } \right|^{2} \right) \le 1\]
is satisfied; note that for such a permutation it should firstly reduce the surface integral into the multiple one, possibly in improper sense, then make a rearrangement of the order of integrations (in this case, since the element of area is independent on $\bar{\alpha }$, and only the region of integration depends on it, then the received after the permutation inner integral may be again transformed into a former surface integral ). Lemma 6 follows  from this relation.

\textbf{Lemma 7.} Let $(x_{1} ,y_{1} ,...,x_{2k} ,y_{2k} )$ be an arbitrary solution of the system (1). Then, for arbitrary real numbers $a$ and $b$ the vector $(x_{1} +a,y_{1} +b,...,x_{2k} +a,y_{2k} +b)$ will also be a solution of the system (1).

\textit{Proof.} Let us take an arbitrary equation of the system (1):

\[x_{1}^{i} y_{1}^{j} +\cdots +x_{k}^{i} y_{k}^{j} -x_{k+1}^{j} y_{k+1}^{j} -\cdots -x_{2k}^{j} y_{2k}^{j} =0.\]
We have:
\begin{equation}\label{12}
\left(x+a\right)^{i} \left(y+b\right)^{j} =\sum _{r=0}^{i}\sum _{t=0}^{j}C_{i}^{r} C_{j}^{t} x^{i-r} y^{j-t} a^{r} b^{t}  .     \end{equation}
Therefore, taking $\varepsilon _{s} =1$ when $s\le k$ and $\varepsilon _{s} =-1$ when $s>k$ we get:
\begin{equation}\label{13}
\sum _{s=1}^{2k}\varepsilon _{s} \left(x_{s} +a\right)^{i} \left(y_{s} +b\right)^{j} =\sum _{r=0}^{i}\sum _{t=0}^{j}C_{i}^{r} C_{j}^{t} \sum _{s=1}^{2k}\varepsilon _{s} x_{s} ^{i-r} y_{s} ^{j-t} a^{r} b^{t}   .    \end{equation}
Inner sum along \textit{s }on the right hand side of (7) is equal to zero for any $r$ and $t$. Lemma 7 is proved.

\textbf{\textit{Lemma 8.}} Let the conditions of Lemma 2 be satisfied, $G_{0} =G_{0} (\bar{x})$ determined by Lemma 1. Then for any vector of a view $\bar{a}=(a,b,...,a,b)$ the following equality is true:

\[G_{0} (\bar{x})=G_{0} (\bar{x}+\bar{a}).\]

\textit{Proof.} From the equation (7) by differentiating with respect to $x_{s} $ and $y_{s} $, we obtain that the rows of the matrix $A_{0} (\bar{x}+\bar{a})$ are a linear combination of the rows of the matrix $A_{0} (\bar{x})$. Then the minors of the matrices$A_{0} (\bar{x}+\bar{a})$ and $A_{0} (\bar{x})$ are the same. Since $G_{0} $ may be represented as a sum of squares of all minors of a maximal order, then the assertion of the Lemma 8 follows from this observation.

\textbf{\textit{Lemma 9}}. Let the conditions of Lemma 1 be satisfied, $\xi _{1}^{0} =\cdots =\xi _{r}^{0} =0$, and the surface $M$\textit{ }determined by the system of equations

\[f_{1} (x_{1} ,...,x_{n} )=0,\]

\[\cdots \]

\[f_{r} (x_{1} ,...,x_{n} )=0,\]
moreover, the functions $f_{j} (\bar{x})$ continuously differentiable in a certain domain $\Omega _{0} $, including  $\Omega $. Let $G=G(\bar{x})$  be the Gram determinant of gradients of functions $f_{j} (\bar{x})$ which is not equal to zero in $\Omega $. Let, further the one-one transformation of coordinates $\bar{x}=\bar{x}'(\bar{\xi })$ maps some domain $\Omega '$ into $\Omega $ with nonsingular Jacoby matrix

\[Q=Q(\bar{\xi })=\left\| \frac{\partial x_{i} }{\partial \xi _{j} } \right\| _{1\le i,j\le n} ,\]
with continuous in  $\Omega $ elements. Then, for any continuous function $f(\bar{x})$ the following formula is true:

\[\int _{M}f(\bar{x}) \frac{ds}{\sqrt{G} } =\int _{M'}\left|\det Q\right|f(\bar{x}(\bar{\xi })) \frac{d\sigma }{\sqrt{G'} } ,\, G'=\det (JQ\cdot {}^{t} Q{}^{t} J),\]
here $M'$ is the prototype of the surface during the transformation, $d\sigma $ is an element of area in the coordinates $\bar{\xi }$, $J$\textit{ }is a Jacoby matrix of the system of functions $f_{j} (\bar{x})$:

\[J=\frac{\partial (f_{1} ,...,f_{r} )}{\partial (x_{1} ,...,x_{n} )} .\]

The proof follows easily from the relation of Lemma 1.

\textit{}
\begin{center}
\textbf{3. Proof of the theorems.}
\end{center}
\textit{}
\textbf{Proof of Theorem 1}. We will use Lemma 4. We conduct the proof in several stages.

1) Let \textit{E } be the set of those points $(x,y)$ of the unite square for which

\[\left(\frac{\partial F}{\partial x} \right)^{2} +\left(\frac{\partial F}{\partial y} \right)^{2} \le \frac{1}{2k} .\]

Denote $\Pi '\left(\bar{\alpha }\right)$, when $\bar{\alpha }$ fixed, a part of the surface  of Lemma 4 lying in the Cartesian product $E^{2k} $. Then, on $\Pi '\left(\bar{\alpha }\right)$ the condition

\[\sum _{j=1}^{2k}\left(\left|\frac{\partial F(x_{j} ,y_{j} )}{\partial x_{j} } \right|^{2} +\left|\frac{\partial F(x_{j} ,y_{j} )}{\partial y_{j} } \right|^{2} \right) \le 1,\]

can be omitted, because for all points $(x_{1} ,y_{1} ,...,x_{2k} ,y_{2k} )\in E^{2k} $ we have:

\[\sum _{j=1}^{2k}\left(\left|\frac{\partial F(x_{j} ,y_{j} )}{\partial x_{j} } \right|^{2} +\left|\frac{\partial F(x_{j} ,y_{j} )}{\partial y_{j} } \right|^{2} \right) \le \sum _{j=1}^{2k}\frac{1}{2k}  \le 1.\]

Therefore,

\[\theta _{k} \ge 2^{N} \pi ^{N/2} \Gamma (N/2+1)\int _{-\infty }^{\infty }\cdots \int _{-\infty }^{\infty }\left(\int _{\Pi '\left(\bar{\alpha }\right)}ds \right)d\alpha _{01} \cdots d\alpha _{nm}   .\]
Here $\Pi '\left(\bar{\alpha }\right)$ is determined by the system (1) only in $E^{2k} $.

2) Let $P\ge 1$ be a natural number. Consider splitting the unite square into small squares with the right upper vertices at the points $(u_{1} ,u_{2} )=(\nu /P,\mu /P)$ when $1\le \nu ,\mu \le P$ are fixed. Rewrite the polynomial $F(x,y)$ in another form using the Taylor's formula at the point $\bar{u}=(u_{1} ,u_{2} ):$
\begin{equation}\label{14}
F(x,y)=\sum _{s_{1} }^{n}\sum _{s_{2} }^{m}\beta \left(\bar{s}\right)(x-u_{1} )^{s_{1} } (y-u_{2} )^{s_{2} } ;\bar{s}=(s_{1} ,s_{2} )
\end{equation}

(14) defines a linear mapping of the space of points $\alpha _{10} ,....,\alpha _{nm} $ into the space of new variables $\beta _{10} ,....,\beta _{nm} ;\beta _{ij} =\beta (i,j).$ Comparing coefficients, we obtain:
\begin{equation}\label{15}
\alpha _{st} =\sum _{p=s}^{n}\sum _{q=t}^{m}(-1)^{p-s+q-t} \beta \left(p,q\right)C_{p}^{s} C_{q}^{t} u_{1}^{p-s} u_{2}^{q-t}  .
\end{equation}
Rearrange pairs $(s,t)$ by this way: if $s+t<s'+t'$ then the pair $(s',t')$ is preceded by a pair $(s,t)$;  if $s+t=s'+t'$ then compare the first and then the second components, and when $s=s',t<t'$ then we will assume that the pair $(s,t)$ precedes the pair $(s',t').$ In such arrangement the relations (14) can be written in the matrix form $\bar{\alpha }=U\bar{\beta }$; the vector $\bar{\alpha }$ has the components $\alpha _{st} $, and the vector $\bar{\beta }$ has the components $\beta \left(s,t\right)$ taken in the just entered descending order. From (14) we conclude that each equation in the right part contains all of the «greater» indices of $(p,q)$ since $(s,t)$, and the coefficient of $\beta (p,q)$ is equal to 1. This means that $U$ is a triangular matrix with zero elements above the main diagonal, and the diagonal elements are equal to 1:

\[U=\left(\begin{array}{cccc} {1} & {0} & {\cdots } & {0} \\ {u_{21} } & {1} & {\cdots } & {0} \\ {\vdots } & {\vdots } & {\ddots } & {0} \\ {u_{N1} } & {u_{N2} } & {\cdots } & {1} \end{array}\right).\]

Therefore, a linear mapping defined above preserves areas.

3) Returning to (14) we obtain:

\[\frac{\partial F}{\partial x} =\sum _{s_{1} }^{n}\sum _{s_{2} }^{m}s_{1} \beta \left(\bar{s}\right)(x-u_{1} )^{s_{1} -1} (y-u_{2} )^{s_{2} } .  \]

Define in the space of variables $\alpha _{10} ,....,\alpha _{nm} $a region $\pi _{\bar{u}} $ as a prototype of the domain of variables $\beta _{10} ,....,\beta _{nm} $ defined by the conditions:

\begin{equation} \label{16} |\beta \left(\bar{s}\right)|\le 0.1cP^{s_{1} +s_{2} -1} ;(s_{1} ,s_{2} )\ne (n,m), \end{equation}
for $s_{1} +s_{2} \ge 1$ and some $c>0$ an exact value of which will set below; for the greatest coefficient we  assume that the inequality

\begin{equation} \label{17} 2^{-1} cP^{n+m-1} \le \beta _{nm} =\alpha _{nm} \le cP^{n+m-1}  \end{equation}
is satisfied. Then,  if $u_{1} -P^{-1} \le x\le u_{1} ,u_{2} -P^{-1} \le y\le u_{2} $ we have:

\[\left(\frac{\partial F}{\partial x} \right)^{2} +\left(\frac{\partial F}{\partial y} \right)^{2} \le \left(\sum _{s_{1} \ge 1}\sum _{s_{2} }s_{1} cP^{s_{1} +s_{2} -1} P^{-(s_{1} +s_{2} -1)}   \right)^{2} +\]
\[\left(\sum _{s_{1} \ge 0}\sum _{s_{2} \ge 1}s_{2} cP^{s_{1} +s_{2} -1} P^{-(s_{1} +s_{2} -1)}   \right)^{2} \le (n^{2} +m^{2} )\left(nmc\right)^{2} .\]

If we set

\[A=\frac{1}{nm\sqrt{2k(n^{2} +m^{2} )} } \]
then under the conditions (10) and (11) the points $(x,y)$ belong to \textit{E}, because in this case,

\[\left(\frac{\partial F}{\partial x} \right)^{2} +\left(\frac{\partial F}{\partial y} \right)^{2} \le \frac{1}{2k} .\]

It follows from the said above that the volume of the region $\pi _{\bar{u}} $ is equal to the volume of the image, i.e., to the expression

\[0.5cP^{n+m-1} \prod _{(i,j)\ne (0,0),(n,m)}(c/10)P^{i+j-1} =\]
\begin{equation} \label{18}
5(c/10)^{N} P^{0.5(n+m)(n+1)(m+1)-N}
\end{equation}

4) Taking $(u_{1} ,u_{2} )=(\nu /P,\mu /P)$ from the equation (8) we obtain

\[\alpha \left(n-1,m\right)=-u_{1} n\beta \left(n,m\right)+\beta \left(n-1,m\right),\]

\[\alpha '\left(n-1,m\right)=-u'_{1} n\beta \left(n,m\right)+\beta '\left(n-1,m\right),\]

\[\alpha \left(n-1,m\right)-\alpha '\left(n-1,m\right)=\left(u'_{1} -u_{1} \right)n\beta \left(n,m\right)+R,\]
where

\[\left|R\right|\le 2\max \left(\left|\beta \left(n-1,m\right)\right|,\left|\beta '\left(n-1,m\right)\right|\right).\]

Let $u'_{1} $ and $u_{1} $ are fractions of a view

\[u_{1} =\frac{\nu _{1} }{P} ,\; u'_{1} =\frac{\nu '_{1} }{P} ,\; \nu _{1} \ne \nu '_{1} ,\]
and $1\le \nu _{1} \le P,1\le \nu '_{1} \le P,\; \nu _{1} ,\nu '_{1} ,P$ are natural numbers. If

\[|\beta \left(\bar{s}\right)|\le 0.1cP^{s_{1} +s_{2} -1} \]
when $\bar{s}=\left(n,m\right)$, then we will have

\[\alpha \left(n-1,m\right)-\alpha '\left(n-1,m\right)=\frac{\nu '_{1} -\nu _{1} }{P} n\beta \left(n,m\right)+\theta \frac{c}{5} P^{n+m-1} ,\]
where $\left|\theta \right|<1,\beta \left(n,m\right)=\alpha \left(n,m\right)=\alpha _{nm} >\frac{c}{2} P^{n+m-1} .$ Since  $\left|\nu _{1} -\nu '_{1} \right|\ge 1$, then from this it follows that $\alpha \left(n-1,m\right)\ne \alpha '\left(n-1,m\right).$ Otherwise, we have, according to said above,

\[\frac{cn}{2} P^{n+m-2} =\frac{\left|\nu '_{1} -\nu _{1} \right|}{P} n\left|\beta \left(n,m\right)\right|\le \frac{c}{5} P^{n+m-1} ,\]
what is impossible. Then in the assumptions made above the decomposition (8) at different points $\bar{u}$ and $\bar{u'}$ corresponds to the different values of $\bar{\alpha }$, i.e. .

5) If the integers $P$ and $P'$ are such that $1\le P'\le 0.5P$ then any two of the regions $\pi _{\bar{u}'} (P')$ and $\pi _{\bar{u}} =\pi _{\bar{u}} (P)$ also have not common points. This follows from the fact that the  regions of variation of  coordinates $\alpha _{nm} =\beta _{nm} $, i.e. the intervals

\[0.5c\left(P'\right)^{n+m-1} <\alpha _{nm} \le c\left(P'\right)^{n+m-1} ,\] \[c\left(P'\right)^{n+m-1} \le 0.5c\left(P'\right)^{n+m-1} <\alpha _{nm} \le c\left(P'\right)^{n+m-1} \]
does not have common points.

6) Let $l=1,2,3,...\; ,P_{l} =2^{l} .$ From 1), 4) and 5) it follows that in the space $R^{N} $ the domains $\pi _{\bar{u}} \left(P_{l} \right)=\pi _{\nu ,\mu } \left(P_{l} \right),1\le \nu ,\mu \le P_{l} $, $l=1,2,...,$ are distinct. So, we have:

\[\theta _{k} \ge 2^{N} \pi ^{N/2} \Gamma (N/2+1)\sum _{l=1}^{\infty }\sum _{\mu =1}^{P_{l} }\sum _{\nu =1}^{P_{l} }\int _{\pi _{\mu ,\nu } (P_{l} )}\left(\int _{\Pi '\left(\bar{\alpha }\right)}ds \right)d\bar{\alpha }    .\]

Let us make a linear change of variables in the multiple integral along $\bar{\alpha }$ $\bar{\alpha }=U\bar{\beta }$ with the Jacobian 1. Let $\pi '_{\nu ,\mu } \left(P_{l} \right),1\le \nu ,\mu \le P_{l} $ be a prototype of $\pi _{\bar{u}} \left(P_{l} \right)$. According to the previous estimation we have:
\begin{equation} \label{19}
\theta _{k} \ge 2^{N} \pi ^{N/2} \Gamma (N/2+1)\sum _{l=1}^{\infty }\sum _{\mu =1}^{P_{l} }\sum _{\nu =1}^{P_{l} }\int _{\pi '_{\mu ,\nu } (P_{l} )}\left(\int _{\Pi '\left(\bar{\alpha }\right)}ds \right)d\bar{\beta }
\end{equation}

7) The domain of variation of  variables $x_{j} ,y_{j} $ on $\Pi '\left(\bar{\alpha }\right)$ (which now designate as $\Pi '\left(\bar{\beta }\right)$) are bounded by the inequalities $\frac{\mu -1}{P_{l} } <x_{j} <\frac{\mu }{P_{l} } ,$ $\frac{\nu -1}{P_{l} } <y_{j} <\frac{\nu }{P_{l} } ,j=1,...,2k$. Then, we get some piece $\Pi '_{l,\mu ,\nu } \left(\bar{\beta }\right)$ of the surface $\Pi '\left(\bar{\beta }\right)$. For each $\bar{\beta }$ we only diminish, according to 3), the internal surface integral in (13), replacing it by the integral taken over this piece. Thus, we obtain:

\[\theta _{k} \ge 2^{N} \pi ^{N/2} \Gamma (N/2+1)\sum _{l=1}^{\infty }\sum _{\mu =1}^{P_{l} }\sum _{\nu =1}^{P_{l} }\int _{\pi '_{\mu ,\nu } (P_{l} )}\left(\int _{\Pi '_{l,\mu ,\nu } \left(\bar{\beta }\right)}ds \right)d\bar{\beta } .\quad    \]

8) By  Lemma 7, to every solution of the system (1) correspond the series of solutions

\[\bar{x}'_{s} =\bar{x}_{s} +\bar{a}(\mu ,\nu ),\; s=1,...,2k,\]
where$\bar{a}(\mu ,\nu )$ is an arbitrary vector of the form

\[\left(-\frac{\mu -1}{P_{l} } ,-\frac{\nu -1}{P_{l} } ,...,-\frac{\mu -1}{P_{l} } ,-\frac{\nu -1}{P_{l} } \right),\]
i.e. in the surface integral on $\Pi '_{l,\mu ,\nu } \left(\bar{\beta }\right)$ one can make the translation by the vector $\bar{a}(\mu ,\nu )$. Then we get a piece $\Pi \left(l\right)$ of the surface where the variables are related only by the conditions (1) and the inequalities

\[0\le x_{j} ,y_{j} \le P_{l}^{-1} ,j=1,...,2k.\]
Therefore,
\begin{equation} \label{20}
\theta _{k} \ge 2^{N} \pi ^{N/2} \Gamma (N/2+1)\sum _{l=1}^{\infty }\sum _{\mu =1}^{P_{l} }\sum _{\nu =1}^{P_{l} }\int _{\pi '_{\mu ,\nu } (P_{l} )}\left(\int _{\Pi \left(l\right)}ds \right)d\bar{\beta }
 \end{equation}
Applying Lemma 1 in the integral on $\Pi \left(l\right)$, let us make a change of variables

\[x_{j} =P_{l}^{-1} z_{j} ,y_{j} =P_{l}^{-1} w_{j} ,j=1,...,2k.\]

The matrix Q in the transformation is a diagonal with the diagonal elements equal $P_{l}^{-1} ,$ i.e., $Q=P_{l}^{-1} \cdot I,$ where $I$ is the identity matrix of order 4\textit{k}. Denote $\Pi _{0} $ the surface determined by the system (1) and the conditions $0\le z_{j} \le 1,\; 0\le w_{j} \le 1,\; j=1,...,2k$ and assume that $f(\bar{x})=G^{1/2} =$ $=G^{1/2} (\bar{x})$. Then  we get, according to Lemma 2:

\[\int _{\Pi \left(l\right)}ds =\int _{\Pi \left(l\right)}f(\bar{x}) \frac{ds}{\sqrt{G} } =P_{l}^{-4k+N} \int _{\Pi _{0} }ds ,\]
so, in the notation of Lemma 2,

\[G'=\det (BQ\cdot {}^{t} Q{}^{t} B)=P_{l}^{-2N} \det (B\cdot {}^{t} B)=P_{l}^{-2N} G;\]
here $B$ denotes the Jacobi matrix of the system of functions, standing on the left parts of the system (1). By $\eta $ denote the integral

\[\eta =\int _{\Pi _{0} }ds. \]

This integral depends of the system (1) only. Now, from (20) and (18) one has

\[\theta _{k} \ge 2^{N} \pi ^{N/2} \Gamma (N/2+1)\sum _{l=1}^{\infty }P_{l}^{-4k+N+2} \int _{\left|\beta _{10} \right|\le \frac{c}{2} P_{l}^{n+m-1} }d\beta _{N} \cdots \int _{\left|\beta _{nm} \right|\le \frac{c}{2} }d\beta _{1}    \ge \]

\[\ge \eta c(N)\sum _{l=1}^{\infty }\left(2^{l} \right)^{-4k+N+2+\sum _{i=0}^{n}\sum _{j=0,i+j>0}^{m}(i+j-1)  }  =\]
\[=c(N)\sum _{l=1}^{\infty }\left(2^{l} \right)^{-4k+2+(n+m)(n+1)(m+1)/2}  ;\]

\[c(N)=5\cdot 2^{N} \pi ^{N/2} \Gamma (1+N/2)\left(10nm\sqrt{2k(n^{2} +m^{2} )} \right)^{-N} .\]

A series along $l$ and , therefore, the integral is divergent when

\[4k\le 2+(n+m)(n+1)(m+1)/2,\]

only if $\eta >0.$

9) We will prove that $\eta >0.$ Consider two cases.

1. If $m+n\ge 4$ then there exists a positive integer $k$ such that

\[1+(n+m)(n+1)(m+1)/4\ge 2k\ge N,\]
because $1-N+(n+m)(n+1)(m+1)/4\ge 2.$ Then, the system (1) has the trivial solution with arbitrary $x_{1} ,y_{1} ,...,x_{k} ,y_{k} $which can be selected so that $G_{0} $ does not vanish. By the theorem on the implicit functions the system (1) defines $4k-N$-dimensional submanifold in some neighborhood of the given point. Therefore, $\eta >0.$

2. Now let $m+n<4.$ Then either $n=m=1$, or $n=2,m=1$ (by the symmetry). We prove the divergence of $\theta _{k} $ in these cases immediately.

The case of  $n=m=1.$ The largest $k$ such that the condition of the theorem are satisfied  is 1. Have

\[F(x,y)=\alpha x+\beta y+\gamma xy.\]

Then,

\[\int _{0}^{1}\int _{0}^{1}e^{2\pi iF(x,y)} dxdy = \int _{0}^{1}\int _{0}^{1}e^{2\pi i\gamma xy} e^{2\pi i(\alpha x+\beta y)} dxdy . \]
Considering the last integral as a Fourier transformation of the functions $e^{2\pi ixy} $ in the unite square, apply the Parseval equality:

\[\int _{-\infty }^{\infty }\int _{-\infty }^{\infty }\left|\int _{0}^{1}\int _{0}^{1}e^{2\pi iF(x,y)} dxdy  \right| ^{2} d\alpha d\beta =4\pi ^{2}  \int _{0}^{1}\int _{0}^{1}\left|e^{2\pi i\gamma xy} \right|^{2} dxdy  =4\pi ^{2} .\]
By integrating with respect to $\gamma $, we obtain the divergent integral.

In the case $n=2,m=1$ we have

\[F(x,y)=\alpha x+\beta y+\gamma xy+\delta x^{2} +\varepsilon x^{2} y.\]
Let us consider the integral

\[\left(\int _{0}^{1}\int _{0}^{1}e^{2\pi iF(x,y)} dxdy  \right)^{2} =\]
\begin{equation} \label{21}
=\int _{0}^{1}\int _{0}^{1}\int _{0}^{1}\int _{0}^{1}e^{2\pi i\left[\alpha (x+u)+\beta (y+v)+\gamma (xy+uv)+\delta (x^{2} y+u^{2} v)\right]} dxdydudv.    \end{equation}

Let's make the change of variables

\[t_{1} =x+u,t_{2} =y+v,t_{3} =xy+uv,t_{4} =x^{2} +y^{2} .\]
Jacoby matrix of this system of functions with respect to the variables $x,y,u,v$ is equal

\[\left(\begin{array}{cccc} {1} & {0} & {1} & {0} \\ {0} & {1} & {0} & {1} \\ {y} & {x} & {v} & {u} \\ {2x} & {0} & {2u} & {0} \end{array}\right).\]
The determinant $D$ of this matrix is not equal to zero, when $u\ne x$. Therefore, the change of variables is one to one  in outside of the hyperplane $u=x$ in the four dimensional space of variables $x,y,u,v$. Thus, [see [17, problem 13.30] is the Fourier transformation of the function $D^{-1} e^{2\pi i(x^{2} y+u^{2} v)} $, in the cube $[0,2]^{4} $, if we consider it as a function of variables $t_{1} ,t_{2} ,t_{3} ,t_{4} $:

\[\left(\int _{0}^{1}\int _{0}^{1}e^{2\pi iF(x,y)} dxdy  \right)^{2} =\]
\[\int _{0}^{2}\int _{0}^{2}\int _{0}^{2}\int _{0}^{2}e^{2\pi i\left[\alpha t_{1} +\beta t_{2} +\gamma t_{3} +\delta t_{4} \right]} D^{-1} e^{2\pi i(x^{2} y+u^{2} v)} dt_{1} dt_{2} dt_{3} dt_{4}     .\]

From here, as above, we conclude that

\[\int _{-\infty }^{\infty }\cdots \int _{-\infty }^{\infty }\left|\int _{0}^{1}\int _{0}^{1}e^{2\pi iF(x,y)} dxdy  \right| ^{4} d\alpha \cdots d\delta =\]
\[=4\pi ^{2}  \int _{0}^{2}\int _{0}^{2}\int _{0}^{2}\int _{0}^{2}\frac{dt_{1} dt_{2} dt_{3} dt_{4} }{D^{2} } >0    .\]
when $\varepsilon \in R$. Integrating with respect to $\varepsilon $, we get a divergent integral. Theorem 1 is proved.

\textbf{Proof of Theorem 2.} Let \textit{k} be an integer so that $4k>2+(n+m)(n+1)(m+1)/2$. By Lemma 1

\[\theta _{k} =\left(2\pi \right)^{N} \int _{\Pi }\frac{ds}{\sqrt{G_{0} } }  ,\]
where the surface integral is taken along the surface $\Pi $, determined by the system of equations (1) in 4\textit{k}-dimensional unite cube. Note, also, that the integral on the right part is understood in improper sense. In consent with the lemma 4 equation $G_{0} =0$ determines in $R^{4k}$ some real algebraic equation. As it was shown in [46] this equation can be solved with respect any of variables and the algebraic set defined by this equation is placed on finite number of surfaces. Consider the system of equations (1). This system of equations defines a submanifold $\Pi _{0} $ in the open domain defined by deleting the algebraic set determined by the equation above. Since we can resolve this equation with respect any of variables then the considered equation $G_{0} =0$ determines on the surface $\Pi _{0} $ a sub variety of a smaller dimension being a closed subset in $\mathbb{R}^{2k}$, and this subset has a zero Jordan measure (see. [23, p. 38]). We pass now to the proof of convergence of the surface integral out of this subset in improper meaning defined in the lemma 2.

1) Let $G=\mathop{\max }\limits_{\bar{x}\in \Omega } G_{0} (\bar{x})$. Let $F(\bar{x})$ be a polynomial of degree $d$. Then, (see. [17 the problem 13.30]) the following inequality shows that $G$ is bounded:

\[G_{0} =\det (A_{0} \cdot {}^{t} A_{0} )\le (4kd^{2} )^{N} .\]

Define subdomains $\Omega _{p} =\{ \bar{x}\in \Omega |2^{-2p} G\le G_{0} $$\le 2^{2-2p} G\} ,p=1,2,...$. Let, further, $\Pi _{p} $ is a part of the surface $\Pi $, lying in $\Omega _{p} $. We take an arbitrary closed region, lying in the domain where $G>0$. In this domain $G_{0} $ reaches its lower bound and, therefore, this domain lies in the union of a finite number of domains $\Omega _{p} $. Then, the convergence of $\theta _{k} $ follows from the convergence of the series $\sum _{p=1}^{\infty }I_{p}  $with $I_{p} =\int _{\Pi _{p} }\frac{ds}{\sqrt{G_{0} } }  $.

2) Estimate the surface integrals $I_{p} $. We perform in the surface integral $I_{p} $ the change of variables $x_{j} =(2h)^{\alpha } u_{j} ,$ $y_{j} =(2h)^{\alpha } v_{j} ;$ $h=\sqrt{G} \cdot 2^{-p} ,$ $\alpha ^{-1} =1+(n+m-2)(n+1)(m+1)/2$, $j=1,...,2k$. Jacobi matrix $Q$ of this transformation is a diagonal matrix of order \textit{4k} the determinant of which is equal to . Since, by Lemma 2, $G'_{0} =\det (A_{0} Q\cdot {}^{t} Q{}^{t} A_{0} )$, then we

have

\[\int _{\Pi _{p} }\frac{ds}{\sqrt{G_{0} } }  =\int _{\Pi '_{p} }\left|\det Q\right|\frac{d\sigma }{\sqrt{G'_{0} } }  ,\]
where $\Pi '_{p} $ is  prototype of the surface $\Pi _{p} $, $d\sigma $ is an element of area on $\Pi '_{p} $. From the view of the matrix $Q$ it follows that $A_{0} Q=\left(2h\right)^{\alpha } A_{0} $and why

\[G'_{0} =\det (A_{0} Q\cdot {}^{t} Q{}^{t} A_{0} )=\left(2h\right)^{2N\alpha } \det (A_{0} \cdot {}^{t} A_{0} ).\]

Each row of the matrix $A_{0} =A_{0} (\bar{x})$ contains monomials with equal degree  and therefore, after of change of varables $x_{j} =(2h)^{\alpha } u_{j} ,$$y_{j} =(2h)^{\alpha } v_{j} $ we get:
\begin{equation} \label{22}
A_{0} =\left(\begin{array}{ccc} {1} &  {\cdots } & {0} \\ {0} &  {\cdots } & {-1} \\ {\cdots } &  {\cdots } & {\cdots } \\ {\left(2h\right)^{\alpha (i+j-1)} iu_{1}^{i-1} v_{1}^{j} } & {\cdots } & {-\left(2h\right)^{\alpha (i+j-1)} ju_{2k}^{i} v_{2k}^{j-1} } \\ {\cdots } & {\ddots } & {\cdots } \\ {\left(2h\right)^{\alpha (n+m-1)} nu_{1}^{n-1} v_{1}^{m} } & {\cdots } & {-\left(2h\right)^{\alpha (n+m-1)} mu_{2k}^{n-1} v_{2k}^{m-1} } \end{array}\right)
\end{equation}

Then, from the equation $\sum _{i+j>0}(i+j-1)=(m+1)\sum _{i=0}^{n}i+(n+1)\sum _{j=0}^{m}j-N=\alpha ^{-1}    $ we get:

\[\det (A_{0} \cdot {}^{t} A_{0} )=\left(2h\right)^{2} G''_{0} ,\]
and $G''_{0} $ has the same view as $G_{0} $, but only the variables $u_{j} ,v_{j} $ vary within $0\le u_{j} ,v{}_{j} \le (2h)^{-\alpha } $. Then,
\begin{equation} \label{23}
I_{p} =\left(2h\right)^{4k\alpha -\alpha N-1} \int _{\Pi '_{p} }\frac{d\sigma }{\sqrt{G''_{0} } } ;\frac{1}{2} \le \sqrt{G''_{0} } \le 1
\end{equation}
Below we denote the matrix (16), after of reducing of all elements of the lines by the common factors $(2h)^{\alpha (i+j-1)} $, again as $A_{0} =A_{0} (\bar{u})$.

3) Fix $p$, and consequently, $h$. Let us denote $\Omega (h)$ the part of $[0,(2h)^{-\alpha } ]^{4k} $ where the inequality $0.5\le \sqrt{G''_{0} } \le 1$ is satisfied. In order to estimate the surface integral on the right hand side of (23), we will dissect, in advance, the surface integral into the parts so that the value of integral can be estimated by using of appropriate projection, comparable in size with pieces of the relevant parts of the surface. These parts are defined by the maximal minors of the  Jacobi matrix $A_{0} $ of the system. Let's dissect $\Omega \left(h\right)$ into no than $t=C_{4k}^{N} $ subdomains $\Omega ^{(\nu )} ,\nu =1,...,t$, overlapping only with their bounds, in each of which one of the minors of the  matrix $A_{0} =A_{0} (\bar{u})$, has a maximal value of modulus among all minors. In each subdomain $\Omega ^{(\nu )} $, at the same time, satisfied the inequality $0.5\le \sqrt{G''_{0} } \le 1$ and the minor with a number $\nu $ accepts everywhere by modulus the maximal values. These subdomains may not be one-connected which hinders one to one continuation of the «local» solutions of the system. The following assertion is a consequence of Lemma 7.

\textit{Each subdomain }$\Omega ^{(\nu )} $ \textit{is a closed set and can be represented as a union of a finite number of one-connected closed domains, as a set of solutions in }$[0,(2h)^{-\alpha } ]^{4k} $\textit{ of the system of polynomial inequalities.}

From said above, it follows that each \textit{subdomain }$\Omega ^{(\nu )} $ is represented as a union $\Omega _{}^{(\nu )} $$=\bigcup _{c\le T_{0} }\Omega (\nu ,c) $, where $\Omega (\nu ,c)$ are one connected subdomains ($T_{0} $ does not depend on $p$). Then, taking among all subdomains that one in which is contained the maximal piece $\pi \left(\nu \right)$ of the surface $\Pi '_{p} $, we can write
\begin{equation} \label{24}
I_{p} \le tT_{0} \left(2h\right)^{4k\alpha -\alpha N-1} \int _{\Pi '_{p} ,0.25\le G''_{0} \le 1}\frac{d\sigma }{\sqrt{G''_{0} } } \le 2tT_{0} (2h)^{4k\alpha -\alpha N-1} \int _{\pi \left(\nu \right)}d\sigma
\end{equation}

4) Let us consider an arbitrary subdomain $\Omega (\nu ,c)$. Although, in it the system allows, in an arbitrary neighborhood of taken solution, the unique solvability with respect to the \textit{same variables,} these solutions in the all subdomain can have several sheets. We estimate the number of sheets of the solvability.

The system (1) admits a unique solvability, say, with respect to the variables $u_{1} ,v_{1} ,....,u_{g} ,v_{g} $ in a neighborhood of a solution $\bar{u}_{0} $ with $g=[N/2]+1$ (if $N$ is odd then $v_{g} $ can be omitted in this list). The remaining variables which we designate as $\xi _{1} ,...,\xi _{4k-N} $, are free, and let $\delta (\bar{u}_{0} )$ be the domain of their variation. Let us denote $\omega \left(\nu \right)$ the union of all domains $\delta (\bar{u}_{0} )$ corresponding all possible points of $\bar{u}_{0} \in \pi \left(\nu \right)$. The mapping $\varphi :\pi \left(\nu \right)\to \omega \left(\nu \right)$ such that $\varphi (\bar{u}_{0} )=\bar{\xi }_{0} $, $\bar{\xi }_{0} =(\xi _{1}^{0} ,...,\xi _{4k-N}^{0} )$ determines, according to Lemma 1 of [24, p. 538], an $f$- sheeted covering. Then, for the system (1) in $\Omega (\nu ,c)$ we have $f$- sheeted solvability. Because, the system (1) is a polynomial, then considering the resultant we by a consecutive elimination of unknowns in the system (1), obtain that $f$ does not exceed some constant $T>0$ depending only on $m,n,N$ and $k$. Consequently, the domain $\Omega (\nu ,c)$ can be dissected into no more than $T$ subdomains $\Delta _{\mu } ,\mu =1,...,f,f\le T$, in each of which the system (1) admits one sheeted solvability.

5) We denote by $\pi _{1} $ the part $\pi _{1} \subset \pi \left(\nu \right)$ that is in the domain $\Delta _{1} $. By Lemma 5, the system (1) together with each solution $\bar{u}=(u_{1} ,...,v_{2k} )$, has also many other solutions of the form $\bar{u}'=(u'_{1} ,...,v'_{2k} )$, with $\bar{u}'=\bar{u}+\bar{a},\; $ $\bar{a}=(a,b,...,a,b),$ $a,b\in {\bf {\rm R}}$. When $\bar{u}=(u_{1} ,...,v_{2k} )\in \pi _{1} $ the  point $\bar{u}'$, with arbitrary real $a,b\in {\bf {\rm R}}$ may lie on $\pi _{1} $, if  $\left|a\right|,\left|b\right|\le (2h)^{-\alpha } $ (at large absolute values of these parameters $\bar{u}'$ no longer belongs to $[0,(2h)^{-\alpha } ]^{4k} $, and, obviously, is not on $\pi $). Further, $G(\bar{u})=G(\bar{u}+\bar{a})$ according to Lemma 6, at $\bar{u}\in \pi _{1} $. The set of all vectors $\bar{a}$, with $a,b\in {\bf {\rm R}}$,  forms a two-dimensional subspace in ${\bf {\rm R}}^{4k} $ which designate as $V$. In $\Omega (\nu ,c)$, and therefore , on $\Delta _{1} $ also, one of the minors of Jacobi matrix, accepts maximal absolute values. Let such a minor will, for example, the minor in which the columns are obtained by differentiating with respect to the first $N$ variable $u_{1} ,v_{1} ,...$. Then, the surface $\pi _{1} $ has a parametric representation:

\[u_{1} =u_{1} \left(\xi _{1} ,\xi _{2} ,...,\xi _{4k-N} \right),\]

\[\cdots \quad \cdots \quad \cdots \]

\[u_{N} =u_{N} \left(\xi _{1} ,\xi _{2} ,...,\xi _{4k-N} \right).\]

Define on $\pi _{1} $ the equivalence relation, counting $\bar{u}\equiv \bar{u}'$, then and only then, when $\bar{u}-\bar{u}'\in V$. Each equivalence class is uniquely determined by an its arbitrary element and consists of all vectors that gotten from this element as results of arbitrarily shiftings by vectors from $V$. The union of all classes of equivalence coincides with the set of solutions (by the property of equivalence relation). It is obvious that each class represents a linear submanifold of a view $\bar{u}+V$, where $\bar{u}$ is an arbitrary solution of the system (1). Conveniently, also, to consider this relation firstly in ${\bf {\rm R}}^{4k-N} $. Then, all the manifold of solutions of the system (1) can be represented as $\pi '_{1} +V$taking the solution $\pi '_{1} $ (it can be identified with a factor set) passing through the fixed point $\bar{u}_{0} $. Indeed, let $\bar{u}$ be an arbitrary solution of system (1). Define the numbers $a=\xi _{4k-N-1} -\xi _{4k-N-1}^{0} $, $b=\xi _{4k-N} -\xi _{4k-N}^{0} $. Vector $\bar{u}-\bar{a}=\bar{u}_{0} $ being a solution of the system belongs to $\pi '_{1} $. Therefore $\bar{u}\in \bar{a}+\pi '_{1} $, and the surface $\pi '_{1} $ can be uniquely determined, if to fix taken in advance two consecutive variables (for example, the  variables $\xi _{4k-N-1} $ and $\xi _{4k-N} $).

Interesting for us solutions are obtained from the variety of solutions by taking pieces  belonging $\pi _{1} $. If, now, we take a submanifold $\pi '_{1} \subset \pi _{1} $, with fixed values $\xi _{4k-N-1} =\xi _{4r-N-1}^{0} ,$ $\xi _{4k-N} =\xi _{4k-N}^{0} $of dimension $4k-N-2$ with a maximal $4k-N-2$-dimensional volume then $\pi _{1} $ will be covered by a variety getting from the submanifold $\pi '_{1} $ by the help of all parallel translations by vectors $\bar{a}\in V$, with $\left|a\right|,\left|b\right|\le (2h)^{-\alpha } $. Therefore, the area of the infinitely small element of the surface $\pi _{1} $ can be represented in the form $\Delta \bar{a}\cdot \Delta _{V} \sigma $, where $\Delta \bar{a}=\Delta a\Delta b$, аnd $\Delta _{V} \sigma $ means the area of the projection of the infinitely small element of the surface $\pi '_{1} $ (area, which is designated as $\Delta \sigma $) into the subspace orthogonal to the subspace $V$. Therefore $\Delta _{V} \sigma \le \Delta \sigma $, and  according to (1) we have

\[I_{p} \le 2tT_{0} (2h)^{4k\alpha -\alpha N-1} \int _{\pi \left(\nu \right)}d\sigma  \le 2tTT_{0} \left(2h\right)^{4k\alpha -\alpha N-1} \times \]
\begin{equation} \label{25}
\times \int _{\left|a\right|,\left|b\right|\le (2h)^{-\alpha } }dadb\int _{\pi '_{1} }d\sigma \le 2tTT_{0} \left(2h\right)^{4k\alpha -\alpha N-2\alpha -1} \int _{\pi '_{1} }d\sigma
\end{equation}
moreover, the constants on the right side are not dependent on $p$.

6) Estimate the last surface integral in (23). Using elements of the matrix $A_{0} =A_{0} (\bar{u})$ we will set a new block matrix $D_{0} $ as follows. Let $l$ be the largest integer such that $Nl\le 4k-N-2$. From the conditions of the theorem, it follows that $4k>2+(N+1)(m+n)/2$. Therefore, if $m+n\ge 4$ such \textit{l }exists. (If $m+n=3$ then $N=5$, and we can take $l=1$). Form $l-1$ submatrices of order $N$ using pairwisely different columns of the matrix $A_{0} $ and replace these submatrices, as blocks, on the diagonal of  the matrix $D_{0} $; take the $l$-th block of a size $N\times (s-(l-1)N)$ set up, as well as above, from the columns of the matrix $A_{0} $, consequently, does not taking previously used columns. Further, add to the getting matrix from below a matrix of a view

\[\left(\begin{array}{cc} {D'} & {D''} \end{array}\right),\]
where $D'$ is the zero matrix of a size $(s-lN)\times N$, and $D''$ is a matrix of an order $s-lN$ of a view

\[D''=\left(\begin{array}{ccc} {(2h)^{\alpha } \xi _{s-lN+1} } & {\cdots } & {0} \\ {\vdots } & {\ddots } & {\vdots } \\ {0} & {\cdots } & {(2h)^{\alpha } \xi _{s} } \end{array}\right).\]

Thus, $D_{0} $ is a block-matrix of the type

\[D_{0} =\left(\begin{array}{cc} {B_{1} } & {B_{2} } \\ {D'} & {D''} \end{array}\right),\]
where $B_{1} $ and $B_{2} $ are the blocks composed by the columns of the matrix $A_{0} $. It is important that all of the independent variables were presented in the entries of the matrix $D_{0} $ (since it was formed by columns of $A_{0} $ then all of the independent variables will be found in monomials of  highest degree, included in $D_{0} $). Note that the block $D''$ plays an auxiliary role complementing the matrix $D_{0} $ to the square matrix of order $s$, which is also important for estimation of the integral (all elements of the matrix $D''$ will vanish in re-differentiation).

We have:

\[\det D_{0} =\det B_{1} \cdot (2h)^{\alpha (s-lN)} \xi _{s-lN+1} \cdots \xi _{s} .\]

Because, minors of the matrix $B_{1} $ are also minors of the  matrix $A_{0} $ then

\[\left|\det B_{1} \cdot {}^{t} B_{1} \right|\le \left|\det A_{0} \cdot {}^{t} A_{0} \right|.\]

By the condition $\sqrt{G''_{0} } =\left|\det A_{0} \cdot {}^{t} A_{0} \right|^{1/2} \le 1$, taking into account the limits of variations of independent variables $\xi _{1} ,...,\xi _{4k-N-2} $, we obtain $\left|\det D_{0} \right|\le 1$. Therefore, replacing the condition of integration in (3), we only increase the surface integral. Let's evaluate now the surface integral on the right hand side of (23) under of such conditions. To this end, replace it with the multiple integral over the domain of independent variables. This domain gotten by projecting $\pi '_{1} $ into the subspace of the independent variables $\xi _{1} ,...,\xi _{4k-N-2} $ (the variables $\xi _{4k-N-1} $ and $\xi _{4k-N} $ are fixed). Then, passing in the surface integral to the independent variables, we obtain the following inequality, using the relation for the element of area in Lemma 1:
\begin{equation} \label{26}
\int _{\pi '_{1} }d\sigma \le \left(C_{4k}^{N+2} \right)^{1/2}  \int _{\tau ,\left|\det D_{0} \right|\le 1}d\xi _{1} \cdots d\xi _{4k-N-2} .
\end{equation}

7) We will use the scheme of work $[$11$]$. For the volume of $\tau $, firstly note that the  trivial bound $\mu \left(\tau \right)\le (2h)^{-s\alpha } $ for it is true. Since the determinant of the matrix $D_{0} $ is a polynomial, then the set where the equality $\det D_{0} =0$ is satisfied has zero Jordan measure. Therefore, the integral on the right hand side of (26) does not exceed (we estimate a more general integral replacing the condition $\left|\det D_{0} \right|\le 1$ by $\left|\det D_{0} \right|\le H$):

\[\int _{\tau ,\left|\det D_{0} \right|\le H}d\xi _{1} \cdots d\xi _{4k-N-2}  \le \sum _{j=1}^{\infty }E_{j}  ,\]

where

\[E_{j} =\int _{2^{-j} H\le \left|\det D_{0} \right|\le 2^{1-j} H}d\xi _{1} \cdots d\xi _{s}  .\]

Let $\rho _{1} =\rho _{1} (\bar{\xi }'),...,\rho _{s} =\rho _{s} (\bar{\xi }')$, where $\bar{\xi }'=(\xi _{1} ,...,\xi _{4k-N-2} )$ are the singular values of matrix $D_{0} $, $\rho _{1} \ge \cdots \ge \rho _{s} $. Then from the inequality

\[2^{-j} H\le \rho _{1} \cdots \rho _{s} \le \rho _{s} \rho _{1}^{s-1} \]

derive

\[\rho _{s} \ge 2^{-j} \rho _{1}^{1-s} H.\]

Assuming $D_{0} =\left(d_{ij} \right)$, by Shur's Teorem [21, p. 212] we have:

\[\rho _{1}^{2} \le \rho _{1}^{2} +\cdots +\rho _{s}^{2} \le \sum _{i,j}d_{ij}^{2}  \le 2s-lN+(N-2)s(2h)^{2\alpha (-m-n+1)} .\]

Therefore,
\begin{equation} \label{27}
\rho _{s} \ge 2^{-j} H\left[2s+(N-2)s(2h)^{-2\alpha (m+n-1)} \right]^{\, (1-s)/2} .
\end{equation}
Let's estimate now $E_{j} ,\, j=1,2,...$.We have:
\[\frac{E_{j} }{2^{1-j} H} \le \int _{2^{-j} H\le \left|\det D_{0} \right|\le 2^{1-j} H}\frac{d\bar{\xi }'}{\left|\det D_{0} \right|}  \le \]
\begin{equation} \label{28}
 \le c_{0} \int _{2^{-j} H\le \left|\det D_{0} \right|\le 2^{1-j} H}d\bar{\xi }'\int _{\left\| D_{0} \bar{\alpha }\right\| \le 1}d\bar{\alpha },
\end{equation}
where $c_{0} =\pi ^{-s/2} \Gamma \left(1+s/2\right)$. Further, from the inequality

\[1\ge \left\| D_{0} \bar{\alpha }\right\| ^{2} =\left({}^{t} D_{0} \cdot D_{0} \bar{\alpha },\bar{\alpha }\right)\ge \rho _{s}^{2} \left\| \bar{\alpha }\right\| ^{2} \ge \rho _{s}^{2} \left|\alpha _{i} \right|^{2} ;\left\| \bar{\alpha }\right\| ^{2} =\sum _{i=1}^{s}\left|\alpha _{i} \right|^{2}  ,\]
for all \textit{i}, according to (27), it follows an estimation

\[\left|\alpha _{i} \right|\le \left(\sum _{i=1}^{s}\left|\alpha _{i} \right|^{2}  \right)^{1/2} \le \rho _{s}^{-1} \le 2^{j} H^{-1} \lambda ;\lambda =\]
\[=\left[2s+(N-2)s(2h)^{-2\alpha (m+n-1)} \right]^{\, (s-1)/2} ,\]
for the variables of integration in the inner integral in (28). Let us introduce into consideration the ball:

\[K=\left\{\bar{\alpha }|\left\| \bar{\alpha }\right\| \le 2^{j} H^{-1} \lambda \right\}.\]
(28) can be rewritten in the form
\begin{equation} \label{29}
\frac{E_{j} }{2^{1-j} H} \le c_{0} \int _{\tau }d\bar{\xi }'\int _{K,\left\| D_{0} \bar{\alpha }\right\| \le 1}d\bar{\alpha }  =c_{0} \int _{\bar{\alpha }\in K}d\bar{\alpha }\int _{\bar{\xi }'\in \tau ,\left\| D_{0} \bar{\alpha }\right\| \le 1}d\bar{\xi }' .
\end{equation}

From the ball \textit{K} remove all stripes $K_{m} $ ($m=1,...,s$), stipulated by the conditions of

\begin{equation} \label{30} \left|\alpha _{m} \right|\le G_{1}^{-1} s^{-1} 2^{j(1-s)} H^{s} \lambda ^{1-s} (2h)^{s\alpha } , \end{equation}
\begin{equation} \label{31}
\left|\alpha _{i} \right|\le 2^{j} H^{-1} \lambda , i\ne m,
\end{equation}
where $G_{1} >0$ will be defined below. Let's designate $K_{0} =\bigcup _{m=1}^{s}K_{m}  $and estimate the measure  of $K_{0} $:
\begin{equation} \label{32}
\mu \left(K_{0} \right)\le c_{0} HG_{1}^{-1} .
\end{equation}
On the right hand side of (23), dissecting a multiple integral into two parts, and defining the first one by a condition $\bar{\alpha }\in K_{0} $, and the second by a condition $\bar{\alpha }\in K\backslash K_{0} $, we estimate the first integral trivially, using the found estimation:
\begin{equation} \label{33}
c_{0} \int _{\bar{\alpha }\in K_{0} }d\bar{\alpha }\int _{\bar{\xi }'\in \tau ,\left\| D_{0} \bar{\alpha }\right\| \le 1}d\bar{\xi }'  \le c_{0} HG_{1}^{-1} .
\end{equation}

8) Let us make a change of variables

\[\bar{\eta }=D_{0} (\bar{\xi }')\bar{\alpha }\]
in the inner integral on the extreme right hand side of (29) for every fixed $\bar{\alpha }\in K\backslash K_{0} $ (the components of the vector $\bar{\xi }'$ are free variables). The Jacoby matrix of change of variables is equal to the inverse of a matrix

\[J=\frac{\partial (\eta _{1} ,...,\eta _{s} )}{\partial (\xi _{1} ,...,\xi _{s} )} =\left(\begin{array}{ccc} {\frac{\partial D_{0} }{\partial \xi _{1} } \bar{\alpha }} & {\cdots } & {\frac{\partial D_{0} }{\partial \xi _{s} } \bar{\alpha }} \end{array}\right),\]
and $\partial ({}^{t} D_{0} )/\partial \xi _{j} $ denotes the matrix collected by differentiating of all the entries of the matrix $D_{0} (\bar{\xi }')$ with respect to the variables $\xi _{j} $. We have:

\[\int _{\bar{\alpha }\in K\backslash K_{0} }d\bar{\alpha }\int _{\bar{\xi }'\in \tau ,\left\| D_{0} \bar{\alpha }\right\| \le 1}d\bar{\xi }'  =\int _{\bar{\alpha }\in K\backslash K_{0} }d\bar{\alpha }\int _{\bar{\xi }'\in \tau ,\left\| \bar{\eta }\right\| \le 1}\left|J\right|^{-1} d\bar{\eta }  .\]

For each $\bar{\eta }$ denote $\tau \left(\bar{\eta }\right)$ the subset in $K\backslash K_{0} $ all points of which are subordinated to the inequality $\left\| D_{0} (\bar{\xi }')\bar{\alpha }\right\| \le 1$. Changing the order of integrations in the last integral, we obtain the inequality:

\[\int _{\bar{\alpha }\in K\backslash K_{0} }d\bar{\alpha }\int _{\bar{\xi }'\in \tau ,\left\| D_{0} \bar{\alpha }\right\| \le 1}d\bar{\xi }'  \le \int _{\left\| \bar{\eta }\right\| \le 1}d\bar{\eta }\int _{\tau \left(\bar{\eta }\right)}\left|J\right|^{-1} d\bar{\alpha }  ,\]
spreading integral to all of these $\bar{\eta }$. Let us consider the matrix

\[J=\frac{\partial (\eta _{1} ,...,\eta _{s} )}{\partial (\xi _{1} ,...,\xi _{s} )} =\left(\begin{array}{ccc} {\frac{\partial D_{0} }{\partial \xi _{1} } \bar{\alpha }} & {\cdots } & {\frac{\partial D_{0} }{\partial \xi _{s} } \bar{\alpha }} \end{array}\right)\]
as a matrix of a linear transformation mapping every vector $\bar{\beta }\in R^{s} $ to the vector $\left(\begin{array}{ccc} {\frac{\partial D_{0} }{\partial \xi _{1} } \bar{\alpha }} & {\cdots } & {\frac{\partial D_{0} }{\partial \xi _{s} } \bar{\alpha }} \end{array}\right)\bar{\beta }$. Obviously, it is linearly also on $\bar{\alpha }$. Therefore, we have bilinear mapping $\Phi :\left(\bar{\alpha },\bar{\beta }\right)\mapsto J\bar{\beta }$. For every fxed $\bar{\alpha }$, that defines a linear mapping the matrix of which as a matrix with polynomial entries nonsingular everywhere on $\tau \left(\bar{\eta }\right)$, with exception of possibly points from  submanifold  of  a zero Jordan measure.

We introduce a matrix $D_{1} $ (as well $A_{1} $) which we get by arranging of the enties of columns of the matrix $D_{0} $ (respectively $A_{0} $), consequently, in a line with subsequently taking of the transposed Jacoby matrix of an obtained system of functions (this matrix has the size $s\times s^{2} $). For each pair $\left(\bar{\alpha },\bar{\beta }\right)\in R^{2s} $ is satisfied an equality $\Phi \left(\bar{\alpha },\bar{\beta }\right)=D_{1} (\bar{\alpha }\otimes \bar{\beta })$, and if ${}^{t} \bar{\alpha }=\left(\alpha _{1} ,...,\alpha _{s} \right)$ and ${}^{t} \bar{\beta }=\left(\beta _{1} ,...,\beta _{s} \right)$ then the symbol ${}^{t} (\bar{\alpha }\otimes \bar{\beta })$ will mean a tensor product  (see. [20, p. 80 or 21, p. 235]). The components of this product ordered lexicographically, and it induces an order in the set of columns of the matrix $D_{1} $. From said above, we obtain:
\begin{equation} \label{34}
\int _{\left\| \bar{\eta }\right\| \le 1}d\bar{\eta }\int _{\tau \left(\bar{\eta }\right)}\left|J\right|^{-1} d\bar{\alpha }  =c_{0} \int _{\left\| \bar{\eta }\right\| \le 1}d\bar{\eta }\int _{\tau \left(\bar{\eta }\right)}d\bar{\alpha }  \int _{\left\| D_{1} (\bar{\alpha }\otimes \bar{\beta })\right\| \le 1}d\bar{\beta } .
\end{equation}

Consider the inner multiple integral on $\bar{\alpha }$ and $\bar{\beta }$:
\begin{equation} \label{35}
\int _{\tau \left(\bar{\eta }\right),\left\| D_{1} (\bar{\alpha }\otimes \bar{\beta })\right\| \le 1}d\bar{\alpha }d\bar{\beta } .
\end{equation}

Let the singular decomposition of the matrix $D_{1} $ is of the form $D_{1} =Q\Sigma T$, where  $Q$ and $T$  are orthogonal matrices of sizes $s$ and $s^{2} $ respectively, and $\Sigma =(\Sigma _{1} |\Sigma _{2} )$ with a diagonal matrix $\Sigma _{1} $ consisting of the singular values of the matrix$D_{1} $, $\Sigma _{2} $ is a  zero matrix (of course columns of $D_{1} $ can be placed in $\Sigma $ in an arbitrary order). Consider in $R^{s^{2} } $ the $2s$ dimensional manifold:

\[\left(\alpha _{1} \beta _{1} ,...,\alpha _{1} \beta _{s} ,...,\alpha _{s} \beta _{1} ,...,\alpha _{s} \beta _{s} \right).\]

 Exclude from the consideration all the hyperplanes of a view $\alpha _{i} =\alpha _{j} $ and  $\beta _{i} =\beta _{j} $ $(i\ne j)$. The union of them has a zero Jordan mesure and by this reason have not effect to the value of the integral..

Let us consider the integral (35) and make change of variables $t_{i} =\alpha _{i} \beta _{i} ,i=1,...,s$. Before applying Lemma 1, we will spend the following reasoning. From the change of variables we find $\alpha _{i} =t_{i} \beta _{i}^{-1} $. Let us denote  conditionally $\bar{\alpha}=\bar{t}\bar{\beta}^{-1}$. We can write

\[D_{1} (\bar{\alpha }\otimes \bar{\beta })=D_{1} \left(\bar{t}\bar{\beta }^{-1} \otimes \bar{\beta }\right).\]

According to (30), for all points of $K \setminus K_{0}$ and arbitrary $i$ it will be fulfilled the inequality

\[G_{1}^{-1} s^{-1} 2^{j(1-s)} H^{s} \lambda ^{1-s} (2h)^{s\alpha } <\left|\alpha _{i} \right|\le 2^{j} H^{-1} \lambda .\]

Below we will impose on $G_{1}(\xi')$ the condition $G_{1}(\xi')>G_{1}$ , where $G_{1}>0$ is the same constant discussed above, does not dependent on $h$ and be defined more precisely later. Therefore,
\begin{equation} \label{36}
\left|t_{i} \right|2^{-j} H\lambda ^{-1} \le \left|\beta _{j_{i} } \right|\le \left|t_{i} \right|\left|\alpha _{i} \right|^{-1} \le \left|t_{i} \right|sG_{1} 2^{j(1-s)} H^{-s} \lambda ^{s-1} (2h)^{-s\alpha},
\end{equation}

for all \textit{i}. Now apply Lemma 1, making the change of variables as it was defined above:

\[\int _{\tau \left(\bar{\eta }\right)}d\bar{\alpha }\int _{\left\| D_{1} \left(\bar{t}\bar{\beta }^{-1} \otimes \bar{\beta }\right)\right\| \le 1}d\beta =  \]
\begin{equation} \label{37}
=\int d\bar{t}\int _{t_{i} =\alpha _{i} \beta _{i} ,\left\| D_{1} \left(\bar{t}\bar{\beta }^{-1} \otimes \bar{\beta }\right)\right\| \le 1}\frac{ds}{\sqrt{\alpha _{1}^{2} +\beta _{1}^{2} } \cdots \sqrt{\alpha _{s}^{2} +\beta _{s}^{2} } }   ;
\end{equation}

We differ the index \textit{s} in the denominator from the varable \textit{s}. Transforming the surface integral into multiple one we get:

\[\int _{\begin{array}{l} {t_{1} =\alpha _{1} \beta _{i_{1} } ,...,t_{s} =\alpha _{s} \beta _{i_{s} } } \\ {\; \left\| D_{1} \left(\bar{t}\bar{\beta }^{-1} \otimes \bar{\beta }\right)\right\| \le 1} \end{array}}\frac{ds}{\sqrt{\alpha _{1}^{2} +\beta _{i_{1} }^{2} } \cdots \sqrt{\alpha _{s}^{2} +\beta _{i_{s} }^{2} } }  =\int \frac{d\beta _{2} \cdots d\beta _{s} }{\beta _{2} \cdots \beta _{s} }  ,\]
moreover, the bounds of variation of variables $\beta_{i}$ are determined by the inequalities (30).

For estimation from below the multiple integral on the right hand side of (37), change the order of integrations:

\[\int \frac{d\beta _{1} \cdots d\beta _{s} }{\beta _{1} \cdots \beta _{s} }  \int _{\left\| D_{1} \left(\bar{t}\bar{\beta }^{-1} \otimes \bar{\beta }\right)\right\| \le 1}dt_{1} \cdots dt_{s}  .\]

The inner integral can be represented in the form of an integral over the surface of the  linear variety $\bar{t}\bar{\beta}^{-1}\otimes \bar{\beta}$ in $R^{s^{2}}$  of  dimension $s$. The element of the area is of the form

\[U\cdot {}^{t} U,\]
where

\[\left(\begin{array}{ccccccccc} {1} & {\beta _{2} \beta _{1}^{-1} } & {\cdots } & {\beta _{s} \beta _{1}^{-1} } & {0} & {0} & {\cdots } & {0} & {\cdots } \\ {0} & {0} & {\cdots } & {0} & {\beta _{1} \beta _{2}^{-1} } & {1} & {\cdots } & {0} & {\cdots } \\ {\cdots } & {\cdots } & {\cdots } & {\cdots } & {\cdots } & {\cdots } & {\cdots } & {\cdots } & {\cdots } \end{array}\right).\]

Therefore,

\[\int _{\left\| D_{1} \left(\bar{t}\bar{\beta }^{-1} \otimes \bar{\beta }\right)\right\| \le 1}dt_{1} \cdots dt_{s}  =\int _{\left\| D_{1} \bar{x}\right\| \le 1}\frac{ds}{\left|U\cdot {}^{t} U\right|}  ,\]
where the surface integral is taken over the  piece of the surface $\bar{t}\bar{\beta}^{-1}\otimes \bar{\beta}$, satisfying the conditions specified under the integral sign. The matrix $U$ contains the identity submatrix which guarantees the inequality $|U\cdot^{t}U|\geq 1$. Transform the linear manifold $\bar{t}\bar{\beta}^{-1}\otimes \bar{\beta}$, acting to it by the matrix $T$ from the singular decomposition of a matrix $D_{1}$. Since $T$ is an orthogonal matrix, then, the value of the integral does not change after of transformation. So, we have:

\[\int _{\left\| D_{1} \left(\bar{t}\bar{\beta }^{-1} \otimes \bar{\beta }\right)\right\| \le 1}dt_{1} \cdots dt_{s}  =\int _{\left\| \Sigma \bar{u}\right\| \le 1}\frac{d\sigma }{\left|U\cdot {}^{t} U\right|}  \le \int _{\left\| \Sigma \bar{u}\right\| \le 1}d\sigma  ,\]

where $d\sigma $ is an element of area  on the variety $\bar{u}=T\left(\bar{t}\bar{\beta }^{-1} \otimes \bar{\beta }\right)$.

Further, we have,

\begin{equation} \label{38} \int _{\left\| \Sigma \bar{u}\right\| \le 1}d\sigma  =c'\sigma _{1}^{-1} \cdots \sigma _{s}^{-1} =c'\det \left(D_{1} \cdot {}^{t} D_{1} \right)^{-1/2} =c'\delta ^{-1}  \end{equation}
(see [13, p. 148]) ($c'$ is a constant). Using the bounds of  variables (30), we find:

\[\int _{\left|t_{i} \right|2^{-j} H\lambda ^{-1} }^{\left|t_{i} \right|sG_{1} 2^{-j(s-1)} H^{1-s} \lambda ^{s-1} \left(2h\right)^{-s\alpha } }\frac{d\beta }{\beta }  =\log \left(sG_{1} 2^{js} H^{-s} \lambda ^{s} \left(2h\right)^{-s\alpha } \right),\]

\begin{equation} \label{39}
\int _{\tau \left(\bar{\eta }\right)}d\bar{\alpha } \int _{\left\| D_{1} (\bar{\alpha }\otimes \bar{\beta })\right\| \le 1}d\bar{\beta } =
\end{equation}
\[\int _{\left\| \bar{\eta }\right\| \le 1}d\bar{\eta }\int _{\tau \left(\bar{\eta }\right)}d\bar{\alpha }  \int _{\left\| D_{1} (\bar{\alpha }\otimes \bar{\beta })\right\| \le 1}d\bar{\beta } \le c_{0} s\cdot \delta _{1}^{-1} \wp _{j}^{s} ,\]
where $\wp_{j}=1+log(sG_{1}2^{js}H*{-1-s}\lambda^{s}(2h)^{-s\alpha})$ . Then, according to the found estimation (33), we obtain the following bound for the integral over $K\setminus K_{0}$ on the right hand side of (23):

\[\int _{\bar{\alpha }\in K\backslash K_{0} }d\bar{\alpha }\int _{\bar{\xi }'\in \tau ,\left\| D_{0} \bar{\alpha }\right\| \le 1}d\bar{\xi }'  \le c_{0}^{2} s\cdot s!H\left(\sum _{j=1}^{\infty }\wp _{j}^{s} 2^{1-j}  \right)\delta _{1}^{-1} .\]

Taking $X=sG_{1}2^{js}H*{-1-s}\lambda^{s}(2h)^{-s\alpha}$ estimate the sum:

\[\sum _{j=1}^{\infty }\wp _{j}^{s} 2^{1-j}  =2\sum _{j=1}^{s}\left[1+sj\log 2+\log X\right]^{s} 2^{-j}  \]

For the estimation of the sum  note that $1+logX>2s^{2}$, if the function

\[exp\{ s\log (1+sj\log 2+\log X)-0.5j\log 2\} \]
monotonically decreasing as a function of $j$. When $1+logX\leq 2s^{2}$ this function has a maximal value $2^{S}s^{2s}$. Therefore,

\[\sum _{j\ge 0}^{}\left[1+sj\log 2+\log X\right]^{s} 2^{-j/2} 2^{-j/2} \le 2(1+s^{2} +\log X)^{s}  .\]

So,
\begin{equation} \label{40}
\int _{\tau ,\left|\det D_{0} \right|\le H}d\xi _{1} \cdots d\xi _{4k-N-2}  \le 2^{s+3} s^{3s} c_{0}^{2} H\delta _{1}^{-1} \wp ^{s} ;\; \wp =1+s^{2} +\log X.
\end{equation}
It is clear that in determining of the matrix $D_{0}$ we can take the columns of the matrix $A_{0}$, among others, containing the elements, the gradients of which form the columns of the maximal, in the above mentioned sense, minor of the  matrix $A_{1}(\bar{\xi'})$. Then,
\begin{equation} \label{41}
\delta _{1}^{2} \ge \left(C_{4kN}^{s} \right)^{-1} G_{1} (\bar{\xi }');\; G_{1} (\bar{\xi }')=\det (A_{1} (\bar{\xi }')\cdot {}^{t} A_{1} (\bar{\xi }')),
\end{equation}
and one assumes that $G_{1}(\bar{\xi'})$. Differentiation is carried out with respect to the  components of the vector $\bar{\xi'}$ which is determined from the system of considered equations, and, therefore, the rows of the matrix $A_{1}(\bar{\xi'})$ are the linear combinations of the rows of the matrix $A_{1}(\bar{u})$ (see the presentation in the beginning of 9)). The matrix $A_{1}(\bar{u})$ is obtained by differentiation with respect to the components of $\bar{u}$ and, therefore, the matrices $A_{1}(\bar{\xi'})$ and $A_{1}(\bar{u})$ have different sizes: $A_{1}(\bar{\xi'})$ has a size $s\times s^{2}$ , and $A_{1}(\bar{u})$ has a size $4k\times s^{2}$. When differentiating with respect to the components of the vector $\bar{\xi'}$, they arising the complex expressions which includes the partial derivatives of the dependent variables with respect to the independent variables of $\bar{\xi'}$. Below we will replace the obtained estimates (40-41) by estimation which includes submatrices of the matrix $A_{1}(\bar{u})$, the receipt of which is not complicated by the difficulties mentioned above (i.e. the differentiation is made only with respect to the  independent variables of the vector $\bar{u}$.

9) The  matrix $A_{1}(\bar{\xi'})$ can be represented in the form $A_{1}(\bar{u})$, where $A_{1}(\bar{u})$ is the matrix introduced above, and $D(\bar{u})$ is a matrix of a view of Lemma 1:

\[D(\bar{u})=\left(\begin{array}{ccccccc} {1} & {0} & {\cdots } & {0} & {\varphi _{11} } & {\cdots } & {\varphi _{1,4k-s} } \\ {0} & {1} & {\cdots } & {0} & {\varphi _{21} } & {\cdots } & {\varphi _{2,4k-s} } \\ {\vdots } & {\vdots } & {\ddots } & {\vdots } & {\vdots } & {\ddots } & {\vdots } \\ {0} & {0} & {\cdots } & {1} & {\varphi _{s1} } & {\cdots } & {\varphi _{s,4k-s} } \end{array}\right)=\left(E_{s} |\Phi \right),\]
where $E_{s} $ a is the identity matrix of order $s$, and the matrix $\Phi $ has the size $s\times(4k-s)$. Therefore, any minor, for example, the minor $M_{1}$, composed of the first $s$ columns of the matrix $A_{1}(\bar{\xi'})$ can be represented in the form

\[\delta _{1} =\left|\begin{array}{cc} {D(\bar{u})\cdot A_{1}^{s} (\bar{u})} & {0} \\ {\Psi } & {E_{4k-s} } \end{array}\right|,\]
moreover, the matrix $A'_{1}(\bar{u})$ is a rectangular matrix composed of the first \textit{s }columns of the matrix $A_{1}(\bar{u})$, $\Psi$ composed of the last $4k-s$ rows of the matrix $A'_{1}(\bar{u})$. Performing elementary transformation over the last lines of the determinant, we find
\begin{equation} \label{42}
\delta _{1} =\left|\begin{array}{cc} {\left[A_{1}^{s} (\bar{u})\right]} & {-\Phi } \\ {\Psi } & {E_{4k-s} } \end{array}\right|;
\end{equation}
here the matrix $\left[A_{1}^{s} (\bar{u})\right]$ is composed of the first $s$ rows of the matrix $A_{1}^{s} (\bar{u})$ so that the two blocks of the first column of the determinant $\delta _{1} $ form the matrix $A_{1}^{s} (\bar{u})$: $A_{1}^{s} (\bar{u})=\left(\begin{array}{c} {\left[A_{1}^{s} (\bar{u})\right]} \\ {\Psi } \end{array}\right)$. Dissect the area $\tau $ into two parts: in the first one the condition  $G_{1} (\bar{\xi }')\ge G_{1} $ is satisfied, and in the remaining part of the $\tau $ we have $G_{1} (\bar{\xi }')\le G_{1} $. Denoting by $\mu _{1} $ and $\mu _{2} $ the areas of the relevant parts of the surface $\pi '_{1} $, for (26) we obtain:
\begin{equation} \label{43}
\int _{\pi '_{1} }d\sigma \le Z(\mu _{1} +\mu _{2} ) ,
\end{equation}
with the $Z$ not dependent on $p$. To estimate $\mu _{1} $, we  use the relations (23) and (37-40). We have:
\begin{equation} \label{44}
\mu _{1} \le 40\cdot 2^{s} s^{3s+1} c_{0}^{2} tMTT_{0} \left(C_{4k}^{N+2} \right)^{1/2} \left(C_{4kN}^{s} \right)^{1/2} HG_{1}^{-1} \wp ^{s} , \end{equation}
with $M$ not dependent on $p$.

10) Evaluation of the value $\mu _{2} $ can be reduced to the estimation like the estimation of the integral (43) was already obtained above. Firstly define out the area where the condition of the form $\eta \le M_{1} (\bar{\xi }')\le 2\eta $ (the corresponding area designate as $\mu '_{2} $) is satisfied, where $M_{1} =M_{1} (\bar{\xi }')$ is a minor containing elements gradients of which form the maximal minor of the matrix $A_{2} (\bar{\xi }')$. Further,

\[(2\eta )^{-1} \mu _{2}^{(1)} \le \int _{\pi '_{1} ,\eta \le \left|M_{1} \right|\le 2\eta }\frac{1}{\delta _{1} } d\sigma  =c_{0} \int _{\pi '_{1} ,\eta \le \left|M_{1} \right|\le 2\eta }d\sigma \int _{\left\| D_{1} \bar{v}\right\| \le 1}d\bar{v}  ,\]
and the integral in the right part is taken on the part of the product $\pi '_{1} \times R^{4k} $ where the imposed conditions on the variables are satified. For every $\bar{u}\in \pi '_{1} $, all of the functions that depend on $\bar{u}$ continue as a constant by parallel translation by the vectors from the space «orthogonal» to the surface $\pi '_{1} $ (in other words, $\pi '_{1} \times R^{4k} $ formed up by parallel translations of the surface $\pi '_{1} $). By changing the order of integrations, we obtain:
\begin{equation} \label{45}
(2\eta )^{-1} \mu _{2}^{(1)} \le 2c_{0} \int _{\left\| D_{1} \bar{\nu }\right\| \le 1}d\bar{\nu } \int _{\pi '_{1} \left(\bar{\nu }\right),\eta \le \left|M_{1} \right|\le 2\eta }d\sigma  ,
\end{equation}
moreover, $\pi '_{1} (v)$ is a projection of the surface $\pi '_{1} $ corresponding $\bar{v}$ after of changing the order of integration. We  will perform in the inner surface integral the change of variables according to the formula: $\bar{\beta }=D_{1} \bar{v}$, applying Lemma 2. We get, in the terms of this Lemma:
\begin{equation} \label{46}
\int _{\pi '_{1} (\bar{v}),\eta \le \left|M_{1} \right|\le 2\eta }\sqrt{G} \frac{d\sigma }{\sqrt{G} }  =\int _{\omega (\eta )}\left|\det Q\right|\sqrt{G} \frac{d\sigma '}{\sqrt{G'} }  ,
\end{equation}
moreover, $\omega (\lambda )$ is a prototype of the surface  in the considered mapping,

\[G'=\det (JQ\cdot {}^{t} Q{}^{t} J);\]
here $Q$ is a Jacoby matrix of the mapping which is equal to the inverse of a matrix

\[\frac{\partial \bar{\beta }}{\partial \bar{u}} =\frac{\partial (\beta _{1} ,...,\beta _{4k} )}{\partial (u_{1} ,...,u_{4k} )} ,\]

And $J=A_{0} ,0.5\le \sqrt{G} \le 1$, according to (23). As it is known, the rows of the matrix $A_{0} $form a subspace $M$, orthogonal to the span $M'$ of the system of rows of the matrix $D(\bar{u})$. So, $R^{4k} =M\oplus M'$, and, consequently, the element of area can be represented in the form $d\bar{w}=d\bar{y}d\bar{z}$. Each vector $\bar{w}\in R^{4k} $ can be represented as a sum of vectors $\bar{y}$ and $\bar{z}$ from the subspaces $M$ and $M'$, and $\left\|\| Q(\bar{y}+\bar{z})\right\| \le \left\| Q\bar{y}\right\| +\left\| Q\bar{z}\right\|\|$. Let $\left\| \bar{y}_{1} ...,\bar{y}_{l} \right\| $ be the basis consisting of the system of rows of $A_{0} $, $\left\| \bar{z}_{1} ...,\bar{z}_{m} \right\| $ is a basis consisting of the system of rows of $D=D(\bar{u})$. Making the change of variables by  formulas $\bar{w}=W\bar{x}$ where $W=\left(\begin{array}{c} {A_{0} } \\ {D} \end{array}\right)$, we get

\[\left|\det Q^{-1} \right|=c_{0} \int _{\left\| Q\bar{w}\right\| \le 1}d\bar{w}=c_{0} \int _{\left\| Q(A_{0} \bar{u}+D\bar{v})\right\| \le 1}\frac{1}{\sqrt{G} \sqrt{D\cdot {}^{t} D} }  d\bar{u}d\bar{v}\ge  \]

\[\ge c_{0} \int _{\left\| QA_{0} \bar{u}\right\| \le 1/2}\frac{1}{\sqrt{G} } d\bar{u} \int _{\left\| QD\bar{v}\right\| \le 1/2}\frac{1}{\sqrt{D\cdot {}^{t} D} } d\bar{v} .\]

Therefore,

\[\left|\det Q^{-1} \right|=c_{0} \int _{\left\| Q\bar{w}\right\| \le 1}d\bar{w}\ge  \]

\[\ge \frac{2^{-2k} \Gamma (1+2k)}{(2\pi )^{2k} \Gamma \left(1+l/2\right)\Gamma \left(1+2k-l/2\right)\sqrt{D\cdot {}^{t} D} \sqrt{G} } \frac{1}{\sqrt{G'} } \frac{1}{\sqrt{\left|{}^{t} D{}^{t} QQD\right|} } .\]

Substituting in (45), we find:

\[\int _{\pi '_{1} (\bar{v}),\eta \le \left|M_{1} \right|\le 2\eta }d\sigma  \le c'\int _{\omega (\eta )}\sqrt{D\cdot {}^{t} D} \sqrt{{}^{t} D{}^{t} QQD} \cdot Gd\sigma ' .\]

The matrix $QD$, at any point $\bar{x}\in \Pi '_{1} $ being a solution of the system, is a Jacoby matrix of the invers transformation, i. e. $(QD)^{-1} $ \textit{coincides with the Jacoby matrix of the change of variables $\bar{\beta }=D_{1} \bar{\nu }$ with respect to the $\bar{\xi }'$}. Then, denoting $C_{0}$ some constant that does not depending on $p$ and $\eta$, we can come back to (31), and then, with the help of reasoning of the clause 8) from (34) arrive to the relation already obtained above in (39), with the replacement of $\delta_{1}$ by $delta_{2}$, $H$ by $\eta$:

\[\mu _{2}^{(1)} \le 4c_{0} \eta \int _{\left\| D_{1} \bar{\nu }\right\| \le 1}d\bar{\nu } \int _{\pi '_{1} \left(\bar{\nu }\right),\eta \le \left|M_{1} \right|\le 2\eta }d\sigma  .\]
This integral is estimated as above, and one counts the relations $gG_{1}^{-1} <<\left|\det M_{1} \right|<<$ $<<\eta ^{-1}$ with positive $g$, depending only on $N$ and $k$. 

Now we note that in accordance with the lemma 5 matrisies $A_{j}(\bar{u})$ have a maximal rank everywhere, with exception for the points of  some subset of zero Jourdan measure. Then replacing $\eta =G_{1} $ by $\eta /2,\eta /4,...$, and summarizing, we find:

\[\mu _{2} \le T_{1} G_{1} G_{2}^{-1} \wp _{1}^{s} ,\]
where $T_{1}$ is a positive constant, depending only on $N$ and $k$. Note that $\wp_{1}$ is an expression similar to $\wp$. Considering $G_{1}=H^{1/2}
G_{2}^{1/2}$, we find the estimation:

\[\int _{\delta }d\sigma \le ZT_{2} \cdot H^{1/2} G_{2}^{1/2} \wp '^{s}  +Z\mu ''_{2} ,\]
where $T_{2}$ is a positive constant depending only on $N$ and $k$ and $\wp'=max(\wp, \wp_{1}$. Continuing, thus, after a few steps, we arrive at the inequality
\begin{equation} \label{47}
\int _{\delta }d\sigma <<G_{d-1}^{-1/(d-1)} \tilde{\wp }^{s}  , \end{equation}
with $d=n+m$ being equal to the degree of the polynomial $F(x,y)$, $\widetilde{\wp}=max(\wp, \wp_{1},$ $...,\wp_{d-2})$, and the constant hidden under the sign of $<<$ depends only on $d$ and $K$. Note that the values $H, G_{1},...,G_{d-1}$ are defined by the equalities

\[H=1,G_{1} =H^{1/2} G_{2}^{1/2} ,G_{2} =H^{1/3} G_{3}^{1/3} ,...,G_{d-1} =H^{1/(d-1)} G_{d-1}^{(d-2)/(d-1)} .\]

12) Now we estimate $G_{d-1} $ from below. The matrix $D_{d-1}$ composed of the partial derivatives of order \textit{d }of the monomials $u_{1}^{n}v_{1}^{m},...,u_{2k}^{n}v_{2k}^{m}$  with respect to the variables of the vector $]bar{u}$. It is easy to note that each row contains at least one non-zero entry of  a view

\[\frac{\partial ^{d} u^{n} v^{m} }{\partial u^{n} \partial v^{m} } =n!m!.\]

In this case, all the other elements of the column are equal to zero. Therefore, $D_{d-1}\cdot ^{t}D_{d-1}$ is a diagonal matrix. Then,

\[G_{d-1} =\det \left[D\cdot (D_{d-1} \cdot {}^{t} D_{d-1} )\cdot {}^{t} D\right]\ge (n!m!)^{s} .\]

According to (40), taking into account the obvious relation $\widetilde{\wp}<<logh^{-1}$, we have:

\[\int _{\delta }d\sigma <<\left(\log h^{-1} \right)^{s}  .\]
From (17), (18) we derive

\[I_{p} <<\left(2h\right)^{4k\alpha -N\alpha -1} \left(\log h^{-1} \right)^{s} .\]

Since $h=\sqrt{G}\cdot 2^{-p}\leq (4k(n+m))^{N}\cdot 2^{-p}$, the series

\[\sum _{p=1}^{\infty }I_{p} <<\sum _{p=1}^{\infty }2^{-p(4k\alpha -N\alpha -2\alpha -1)} p^{s}   ,\]
and, together with him, the special integral are converging.

Note that for any positive $H$ it takes place the inequality

\[\int _{\delta }d\sigma  <<H^{1/(d-1)} G_{d-1}^{(d-2)/(d-1)} \left(\log h^{-1} \right)^{\, s} ,\]
similar to (31). Therefore, a subset of the surface defined above by the condition $G_{0}=0$, has a zero Lebesgue measure, that we used above. The proof of the theorem 2 is completed.

\textbf{References}

1. Vinogradov I. M. Method of trigonometric sums in the number theory. M., Nauka, 1971 (rus).

2. Vinogradov I. M., Karatsuba A. A. Method of trigonometric sums in the number theory Transactions of  Mathematical Institute of the USSR, v.168, (1984), p. 4-30.

3. Arkhipov G. I., Karatsuba A. A., Chubarikov V. N. Trigonometric integrals, Izv. Academy of Sciences. of  USSR, math. ser. (1979), v.43, №5, pp.971-1003.

4. Arkhipov G. I., Karatsuba A. A., Chubarikov V. N. Theory of multiple trigonometric sums. M., Nauka, 1967.

5. Chubarikov V. N. On multiple trigonometric integrals. Dokl. Academy of Sciences of USSR, (1976), v.227, №6, pp.1308-1310.

6. Chubarikov V. N. On multiple rational trigonometric sums and multiple integrals, Mat. Notes, v.20, №1, (1976), pp. 61-68.

7. Hua Loo Keng On the number of solutions of Tarry's problem Acta Sci. Sinica, (1952), v.1, №1, pp. 1-76.

8. Arkhipov G. I., Karatsuba A. A., Chubarikov V. N. Multiple trigonometric sums and their applications, Proc. MIAN, (1980), v.151, pp.1-128.

9. Arkhipov G. I., Karatsuba A. A., Chubarikov V. N. Theory of multiple trigonometric sums M., Nauka, 1987.

10. Jabbarov I. Sh. On an identity of harmonic analysis and its applications. Dokl. Academy of Sciences of  USSR, (1990), v.314, №5, pp. 1052-1054.

11. Jabbarov I. Sh. On estimates of trigonometric integrals» Transaction of RAS, (1994), v.207, pp. 82-92.

12. Jabbarov I. Sh. On estimates of trigonometric integrals Сhebishevskii sbornik, v. 11, issue 1(33), (2010), pp. 85-108.

13. Bellman R. Introduction to the theory of matrices. M.: Nauka, 1976.

14. Gantmacher F. R. Matrix theory M., Nauka, 1954.

15. Courant R. Course of differential and integral calculus M., Nauka,v.2,1970.

16. Voevodin I. I., Kuznetsov Y. A. Matrices and computing M., Nauka, 1984.

17. Shilov G. E Mathematical analysis. Functions of several real variables M., Nauka, 1972.

18. Bourbaki N. The algebra. Algebraic structures. Linear and polylinear algebra, II-III,  M., Nauka, 1962.

19. Hodge B., Pido D. Methods of algebraic geometry v.2. IL, 1954.

20. Voevodin V. V. Linear algebra, M., Nauka, 1974.

21. Lancaster P. Matrix theory M., Nauka, 1978.

22. Nikolsky S. M. Course of mathematical analysis M., Nauka, v. 1, 4 ed., 1990.

23. Nikolsky S. M. Course of mathematical analysis M., Nauka, v. 2, 4 ed., 1991.

24. Dubrovin B. A., Novikov V. P., Fomenko A. T. Modern geometry. M., Nauka, 2ed., 1968.

25. Evgrafov M. A. Analytic functions. M., Nauka, 1991.

26. Arkhipov G. I. Sadovnichiy V. A. Chubarikov V.N. Lectures on mathematical analysis. M.: Visshay Shcola, 1999.

27. Gelfand A. M. Lectures on linear algebra. M.: Nauka, 1971.

28. Maltsev A. A. Topics in Linear Algebra. M.: Nauka, 1979.

29. Linnik Yu. V. Weyl's sums. Mat. Coll.,1943, v.12, issue 1, pp. 28 - 39.

30. Karatsuba A. A. Mean values of the modulus of trigonometric sums. Izv. AN SSSR, ser. Mat., 37(1973), pp.1203-1227.

31. Ikromov I. A. On the convergence exponent of Trigonometrically integrals Tr. MIAN, (1997), v.218, pp.179-189.

32. Chahkiev M. A. On the exponent of convergence of the singular integral of multidimensional analogue of Tarry's problem Izv. RAN, v. 97, № 2, (2003). pp. 211-224.

33. Landau E. Introduction to differential and integral calculus. M: GIIL, 1948.

34. Danford N., Schwartz J. T. Linear operators. M.: IL, 1962.

35. Hua Loo Keng. Method of trigonometric sums and their application in the theory of numbers,. M: Nauka, 1964.

36. Sprindzhuk V. G. Metric Theory of Diophantine approximation, M. Nauka, 1977.

37. Arkhipov G. I., Karatsuba A. A., Chubarikov V. N. Convergence exponent of special integral of Tarry's problem, Academy of Sciences of USSR, Math. ser. (1979), v.43, №5, pp.971-1003.

38. Chubarikov V. N. On asymptotic formulas for the integral I. M. Vinogradov and its generalizations. Tr. MIAN.(1981) v.157., pp.214-232.

39. Karatsuba A. A. Basics of analytic number theory. M.,Nauka, 1983.

40. Rudin U. Basics of mathematical analysis. M. Mir, 1976.

41. Titchmarsh E. C. Theory of functions. M. GITTL, 1953.

42. Hinchin A. Ya. Short course of mathematical analysis. M.-GITTL, 1953.

43. Fihtenholz G. M. Course of differential and integral calculus. v. 2. M. Fizmatlit, 2006.

44. Arnold V. I., Varchnko A. N., Huseyn-zadeh S. M. Singularities of differentiable mappings. M., Nauka, 2009.

45. Aslanova N. Sh, Jabbarov I. Sh. On an extremality of some algebraic manifolds. Materials of the international conference devoted to 870 year anniversary of the birth of the great poet and philosopher Nizami Ganjavi «Mathematical theories, the problems of their application and teaching», 23-25 September, 2011. Ganja, 2011, pp.103-104.

46. Jabbarov I. Sh. On the structure of some algebraic varieties. Transactions of NAS of Azerbaijan, Issue Mathematics, 36 (1), 74-82 (2016). Series of Physical-Technical and Mathematical Sciences.

Haydar Aliev avenue, 187,

Ganja Stat University,

Ganja, Azerbaijan

\end{document}